\documentclass[12pt, twoside, a4paper]{article}

\usepackage[english]{babel}  
\usepackage{float}
\usepackage[utf8]{inputenc}  
\usepackage[T2A]{fontenc}      
\usepackage{amsmath, amssymb, amsthm, amsfonts, amsthm, dsfont,mathrsfs,wrapfig,graphicx}
\usepackage[]{hyperref}
\usepackage{epigraph}
\usepackage{accents}
\usepackage[LGR,T1]{fontenc} 

\DeclareSymbolFont{upgreek}{LGR}{cmr}{m}{n}
\SetSymbolFont{upgreek}{bold}{LGR}{cmr}{bx}{n}
\DeclareMathSymbol{\Eta}{\mathord}{upgreek}{`H}

\usepackage{abstract}
\usepackage{mathrsfs}
\usepackage[cal=cm,scr=dutchcal]{mathalfa} 
\usepackage{bm}
\DeclareMathAlphabet{\mathpzc}{OT1}{pzc}{m}{it}

\DeclareMathOperator{\supp}{supp}
\DeclareMathOperator{\dist}{dist}

\DeclareMathOperator{\const}{const}
\DeclareMathOperator{\clos}{clos}
\DeclareMathOperator{\loc}{{loc}}

\DeclareMathOperator{\Op}{Op}

\DeclareMathOperator{\Isom}{Isom}
\DeclareMathOperator{\arccosh}{arccosh}

\DeclareMathOperator{\mmod}{mod}
\DeclareMathOperator{\geod}{geod}

\newcommand{\R}{{\mathbb R}}
\newcommand{\HH}{{\mathbb H}}

\newcommand{\argmax}[1]{\underset{#1}{\arg\max}}

\newtheorem{theorem}{Theorem}[section]
\newtheorem{lemma}[theorem]{Lemma}
\newtheorem{define}[theorem]{Definition}
\newtheorem{predl}[theorem]{Proposition}
\newtheorem{sled}[theorem]{Corollary}

\newcommand{\CC}{\mathbb C}

\usepackage{enumitem}

\usepackage{fancyhdr}
\renewenvironment{abstract}
{\small
	\begin{center}
		\bfseries \abstractname\vspace{-.5em}\vspace{0pt}
	\end{center}
	\list{}{%
		\setlength{\leftmargin}{2cm}
		\setlength{\rightmargin}{\leftmargin}%
	}%
	\item\relax}
{\endlist}

\makeatletter
\def\blindfootnote{\gdef\@thefnmark{}\@footnotetext}
\makeatother

\usepackage{mathtools}

\begin{document}

	\LARGE 
	\Large
	\normalsize
	
	\thispagestyle{empty}
	\title{\vspace{-1cm}Growth and nodal current of complexified horocycle eigenfunctions}
	\author{Mikhail Dubashinskiy} 
	\date{\today}

	\newpage
	{\vspace{-6cm}\maketitle}

	\blindfootnote{Chebyshev Laboratory, St.~Petersburg State University, 14th Line 29b, Vasilyevsky Island, Saint~Petersburg 199178, Russia.}
	\blindfootnote{\hspace{0.2mm}e-mail: \href{mailto://mikhail.dubashinskiy@gmail.com}{\texttt{mikhail.dubashinskiy@gmail.com}}}
	\blindfootnote{Research is supported by the Russian Science Foundation grant 19-71-30002.}
	
	\blindfootnote{{Keywords:} quantum ergodicity, semiclassical limit, Fourier integral operator, nodal distribution, horocycle flow, gauge invariance, Grauert tube. 
	}

	\blindfootnote{\hspace{0.0mm}{MSC 2020 Primary}: 58J51; Secondary:  11F37, 37D40, 58J40, 32A60.
	}

	\blindfootnote{This manuscript version is made available under the \href{http://creativecommons.org/licenses/bync-nd/4.0/}{CC-BY-NC-ND 4.0 license}.}

	
	\thispagestyle{empty}

	\renewcommand{\abstractname}{}
	
	\begin{abstract}\vspace{-1.2cm}\footnotesize
		\noindent {\bf Abstract.} 
		We study horocycle eigenfunctions at Lobachevsky plane. They are functions $u\colon \HH=\mathbb C^+=\{z\in\mathbb C\colon \Im z>0\}\to\mathbb C$ such that $\left(-y^2\left(\frac{\partial^2}{\partial x^2}+\frac{\partial^2}{\partial y^2}\right)+ 2i\tau y\frac{\partial}{\partial x}\right)u(x+iy)=s^2 u(x+iy)$, $x+iy\in\mathbb C^+$, with $\tau,s\in\R$,  $\tau$ large and $s/\tau$ small.  In other words, we study eigenfunctions of \emph{magnetic} quantum Hamiltonian on hyperbolic plane.   By Bohr semiclassical correspondence principle, the asymptotic behavior of such functions is related to horocycle flow on $T\HH$. Let $u^\CC$ be analytic continuation of function $u$ to \emph{Grauert tube}; the latter is an open neighbourhood  of $\HH$ in 
		the complexified Lobachevsky plane $\HH^\CC$.  If a sequence of horocycle functions possesses microlocal quantum ergodicity at the admissible energy level (with $\hbar=1/\tau$) then we may find asymptotic distribution of divisor of $u^\CC$. This is done by establishing the asymptotic estimates on $|u^\CC|$ in $\HH^\CC$. Under imaginary-time horocycle flow, microlocalization of $u$ in $T^*\HH$ is taken to localization of $u^\CC$ on $\HH^\CC$. The growth of functions $u^\CC$ as $\tau\to\infty$ turns to be governed by the growth of complexified \emph{gauge factor} occurring in $\tau$-automorphic kernels for functions on $\HH$.
\end{abstract}

{\footnotesize
	\tableofcontents
}

\pagestyle{fancy}

\section{Introduction}

Let $\tau\in\mathbb R$. In hyperbolic Lobachevsky plane $\HH$ implemented as upper half-plane $\CC^+$ consider differential operator  $$D^\tau := -\Delta_\HH+2i\tau y\dfrac{\partial}{\partial x},$$ here $\Delta_\HH:=y^2\left(\dfrac{\partial^2}{\partial x^2}+\dfrac{\partial^2}{\partial y^2}\right)$ is the hyperbolic Laplacian, $x+iy\in\mathbb C^+$. We study asymptotic properties of solutions of eigenfunction equation
$D^\tau u = s^2 u$, $u\colon \HH\to\mathbb C$, for $\tau$ large and $s/\tau$ small.

Let $\tau_n, s_n\in\mathbb R$ 
($n=1,2,\dots$).  Suppose that functions $u_n\colon \HH\to \CC$ are such that 
\begin{equation*}
D^{\tau_n} u_n = s_n^2 u_n
\end{equation*}
and also $\tau_n, 
\tau_n/s_n \xrightarrow{n\to\infty} \infty$. We mostly drop subscript $n$ in what follows.

If we take Planck constant $\hbar=1/\tau$ then the principal symbol of $\tau^{-2}D^\tau$ is  $2H_1-1$ 
where, for $\mathrm b\in \R$, we define \emph{magnetic Hamiltonian}  
$$
H_\mathrm b(x,y,\xi_1,\xi_2) := \frac{(y\xi_1-\mathrm b)^2+(y\xi_2)^2}{2}\colon T^*\HH\to \R
$$ ($z=x+iy\in\HH$, $(\xi_1,\xi_2)$ are cotangent coordinates conjugate to $(x,y)$). Thus, local frequencies of function $u$ with $D^\tau u = s^2 u$ and $s/\tau$ small have to concentrate, as $\tau\to\infty$, near \emph{null level set} $\{H_1=1/2\}\subset T^*\HH$ of the symbol. Notice that, on this set, $H_1$ understood as classical Hamiltonian, generates right \emph{horocycle} flow. 
If we fold $\HH$ into a compact hyperbolic surface by means of an action of a discrete group of isometries then horocycle flow is known to have unique ergodicity property --- unlike geodesic flow which is only ergodic (but is of hyperbolic Anosov type instead). Bohr semiclassical correspondence principle then leads to different conclusions on quantizations of these flows and the stationary states of quantized Hamiltonians.

\begin{define}
	\label{def:QE}
	We say that $\{u_n\}_{n=1}^\infty$ is \emph{quantum ergodic} 
	sequence \emph{(}quantized with $\hbar=1/\tau_n$\emph{)} if, for any $a\in C_0^\infty(T^*\HH)$ understood as a symbol of order $-\infty$,
we have 
\begin{equation*}
	\langle(\Op_{1/\tau_n}a)u_n, u_n\rangle_{L^2(\HH)} \xrightarrow{n\to\infty} \int\limits_{\{H_1=1/2\}} a\, d\mu_L.
\end{equation*}
Here, $\mu_L$ is \emph{horocycle Liouville measure} supported by $\{H_1=1/2\}$, see Section \ref{sec:coords_flows}. \emph{Pseudodifferential operator (PDO)} 
$\Op_{1/\tau_n} a\colon L^2_{\loc}(\HH)\to L^2_{\loc}(\HH)$  is any of semiclassical quantizations of classical observable $a$, see, e.g., \emph{\cite{Zw}} for details. 
\end{define}

Quantum ergodicity of functions $u_n$ means that their local frequencies scaled $1/\tau_n$ times become uniformly distributed at the admissible energy level $\{H_1=1/2\}$. This, of course, depends on the choice of $\hbar_n=1/\tau_n$ in $\Op_{\hbar_n}a$ in  Definition \ref{def:QE} but, in what follows, we do not specify this choice since we always take $\hbar_n=1/\tau_n$. There is a plenty of quantum ergodic sequences, see discussion below in this Introduction.

In the case of free particle on a negatively curved manifold we deal with usual Beltrami--Laplace operator and geodesic flow. Here, the questions on quantum ergodicity for the whole sequence of eignfunctions (\emph{quantum unique ergodicity}) are rather difficult. 
See \cite{Lin06}, \cite{An08}, \cite{DJ}.

Frequency equidistribution of functions $u_n$ leads to consequences on their complexifications. Any real-analytic manifold admits a complexification which is not unique. Any two of such complexifications are biholomorphically equivalent near the original manifold. This is known as Bruhat--Whitney Theorem, see  \cite{BW59}. So, to begin, we may just take $\HH^\CC:=\CC\times\CC$ as complexified hyperbolic plane $\HH$. This set is endowed with \emph{Euclidean} coordinates $(\Re X, \Im X, \Re Y, \Im Y)$, $(X,Y)\in\CC\times\CC$.

Next, we shrink the domain of our interest by replacing $\HH^\CC$ with its  open subset. 
In Section \ref{sec:coords_flows} we define \emph{complex horocycle parametrization mapping}
\begin{equation*}
\R\times\R^+\times (0,1)\times\left(\R\mmod2\pi\right)\ni(x,y,t,\theta)\mapsto  h_{-it}(x+iy,\theta)\in \HH^\CC.
\end{equation*}
Parameter $\theta$ here is 
the slope of a {horocycle} starting from a point $x+iy\in\HH$ and, further, evaluated at imaginary time $-it$ by analytic (with respect to time) continuation; slope $\theta$ is calculated with respect to horizontal line  and at the point $x+iy$. 
Mapping $(x,y,t,\theta)\mapsto h_{-it}(x+iy,\theta)$ with domain as above is injective onto a set of the form $\mathcal G_1\setminus \HH$, $\mathcal G_1\subset \HH^\CC$ being an open vicinity of $\HH$ in $\HH^\CC$ (Proposition~\ref{predl:Grauert}). This set $\mathcal G_1$ is called \emph{radius $1$ horocycle Grauert tube}.  Bijective mapping 
$$
\R\times\R^+\times (0,1)\times\left(\R\mmod2\pi\right)\ni(x,y,t,\theta)\mapsto h_{-it}(x+iy,\theta)\in\mathcal G_1\setminus\HH
$$
gives \emph{horocycle coordinates} $(x,y,t,\theta)$ for punctured Grauert tube  $\mathcal G_1\setminus\HH$. Sometimes we write $t(P)$ and $\theta(P)$ for the latter two coordinates of $P\in\mathcal G_1\setminus\HH$. 

It turns out that such horocycle Grauert tubes, in fact, coincide to the usual geodesic Grauert tubes with recalculated radius (Proposition \ref{predl:Grauert_hor_geod}); this is because the group of isomerties of $\HH$ acts transitively on the spherical bundle $\mathscr S_1\HH$ over $\HH$.

It is easy to see that functions $u_n$ possess analytic continuations $u_n^\CC\colon \mathcal G_1\to \CC$  (Lemma \ref{lemma:SelbergTransform}).  Our first result is on the growth of these complexifications. In horocycle coordinates, define a function $B_0=B_0(x,y,t,\theta)\colon \mathcal G_1\to \R$ as 
\begin{equation}
\label{eq:B0_define}
B_0 := \log\left(\frac{2+(t^2-2t)\cdot(1+\cos\theta)}{2+(t^2+2t)\cdot(1+\cos\theta)}\right) \, \mbox{ on }\, \mathcal G_1\setminus \HH, ~~~ \left.B_0\right|_{\HH} := 0.
\end{equation}
This function is responsible for the growth of $u$ in the following sense:

\begin{theorem}
	\label{th:growth}
	Suppose that $D^{\tau_n}u_n = s_n^2 u_n$ with $\tau_n, s_n\in\R^+$,  $s_n/\tau_n\xrightarrow{n\to\infty}0$ and $\tau_n\xrightarrow{n\to\infty}+\infty$. Assume also that $\sup\limits_{n\in\mathbb N, \, z\in \HH} \|u_n\|_{L^1(\mathcal B_\HH(z,1))}<+\infty$. Here, $\mathcal B_\HH(z,r)\subset\HH$ is the open ball in hyperbolic metric centered in $z\in\HH$ and having radius $r>0$.
	
	Suppose also that sequence $\{u_n\}_{n=1}^\infty$ is quantum ergodic   in the sense of Definition~\ref{def:QE}.
	
	Under these conditions, we have 
	$$
	|\tau_n|^{1/2} \cdot |u_n^\CC|^2 \cdot 
		\exp(|\tau_n| B_0) \xrightharpoondown[\tau\to+\infty]{}^*  b \mbox{ in } \mathcal D'(\mathcal G_1\setminus \HH).
	$$
	Here, $b$  is a smooth function separated from zero on compacts in  $\mathcal G_1\setminus \HH$ and not depending on $\{u_n\}_{n=1}^\infty$. Both sides of the limit relation are understood as densities of measures in 
	$\CC\times\CC\mbox{-Euclidean}$ coordinates in~$\HH^\CC$, and weak* convergence is understood in the sense of distributions.  
\end{theorem}

\noindent {\bf Remark.} Let $\tau<0$. If $D^{\tau}u = s^2 u$ then, for complex conjugates, we have $D^{-\tau}\bar u = \bar s^2 \bar u$. In Corollary \ref{sled:involution} we see that the natural involution 
$$
\HH^\CC\ni(X,Y)\mapsto\imath(X,Y):=(\bar X, \bar Y)\in\HH^\CC
$$ 
preserves $\mathcal G_1$. Thus analytic continuation of $\bar u$ is $\mathcal G_1\ni P\mapsto \bar u(\imath P)$. Therefore, if, in Theorem \ref{th:growth}, we have $\tau_n\to-\infty$ then the conclusion remains true with $B_0$ replaced by $B_0\circ\imath$. The same concerns Theorem \ref{th:main} below.

In what follows we assume that $\tau\ge 0$.

\medskip

\noindent {\bf Remark.} 
In Theorem \ref{th:growth}, we have to cut $\HH$ from $\mathcal G_1$. Lots of our estimates fail when $t$ approaches $0$ (as well as when $t$ is close to $1$). In fact, weak* convergence from Theorem~\ref{th:growth} is valid on each slice 
\begin{equation}
\label{eq:slice_def}
\Sigma_t := \{h_{-it}(x+iy,\theta)\colon x+iy\in\HH,\, \theta\in\R\mmod2\pi\}
\end{equation}
(see Proposition \ref{predl:sliced_growth}) but, as $t\to 0$, this slice tends to $\HH$ and thus degenerates. 
\medskip

Later on in this Introduction, we will discuss the role and meaning of the function $B_0$ giving the answer in Theorem \ref{th:growth}.

\medskip

Now, consider nodal set $\tilde {\mathcal Z}_n := \{P\in\mathcal G_{1}\colon u_n^\CC(P)=0\}\subset \HH^\CC$. Some singularities are possible at this set, but they are always negligible. In all its \emph{non-singular} points set $\tilde {\mathcal Z}_n$ is an analytic submanifold of complex dimension $1$ and thus is canonically endowed with orientation. For any non-singular point  $P\in \tilde Z_{n}$ there is an integer multiplicity of zero of $u_n^\CC$ at $P$, denote this multiplicity by $m_n(P)$. Therefore, $m_n$ and  $\tilde {\mathcal Z}_n$ naturally give rise to de Rham current $\mathcal Z_n$ of dimension $2$: $\mathcal Z_n(\omega) := \int_{\tilde{\mathcal  Z}_n}m_n \omega$ for smooth compactly supported $2$-form $\omega$ in $\mathcal G_{1}$,  $\mathcal Z_n(\cdot)$ denotes application of current $\mathcal Z_n$ to a test form. This current is known to be well-defined.

In a more analytic way, nodal current given by $u_n^\CC$ is equal to the de Rham current defined as 
\begin{equation*}
\mathcal Z_n(\omega) = \dfrac{i}{\pi}\int\limits_{\mathcal G_{1}}\partial\bar{\partial} \log|u_n^\CC|\wedge\omega
\end{equation*}
for test form $\omega$ in $\mathcal G_{1}$. This is known as \emph{Lelong--Poincar\'e formula}. Function $\log|u_n^\CC|$ is understood as $4$-current therein. Operators $\partial$, $\bar{\partial}$ on currents are permanent to those on forms and are given by the complex structure in $\HH^\CC$.  See more in \cite{Chirka}, \cite{LelongGruman}.

In Section \ref{sec:psh}, we take logarithm of the asymptotic relation from  Theorem \ref{th:growth} and derive our second result:

\begin{theorem}
	\label{th:main}
	In the assumptions of Theorem \ref{th:growth}, 
	for nodal currents given by functions $u_n^\CC$, we have 
	$$
	\dfrac{\mathcal Z_{n}}{|\tau_n|} \xrightharpoondown{\tau_n\to+\infty} 
	\frac{i}{2\pi}\bar{\partial}\partial B_0 
	$$
	as de Rham currents of dimension $2$ in $\mathcal G_1$ \emph{(}right-hand side which is a form is also a current\emph{)}.  If $\tau_n\to-\infty$ then $B_0$ is again replaced by $B_0\circ\imath$.
\end{theorem}

Our main example of horocycle quantum ergodic sequence is as follows. Denote by $\Isom^+(\HH)$ the group of orientation-preserving isometries of hyperbolic Lobachevsky plane $\HH$. If $\HH$ is implemented as upper complex half-plane $\mathbb C^+$ then any $\gamma\in \Isom^+(\HH)$ can be written in the \emph{canonical form} $\HH\ni z\mapsto\gamma z = \dfrac{\mathscr az+\mathscr b}{\mathscr cz+\mathscr d}$ for $\mathscr a,\mathscr b,\mathscr c,\mathscr d$ real with $\mathscr a\mathscr d-\mathscr b\mathscr c=1$. Let $\Gamma$ be a discrete torsion-free subgroup in $\Isom^+(\HH)$.  A function $u\colon \HH\to \mathbb C$ is called \emph{$\tau$-form} with respect to $\Gamma$ ($\tau\in\R$) if $u(\gamma z) = \left(\dfrac{\mathscr cz+\mathscr d}{\mathscr c\bar z+\mathscr d}\right)^\tau  u(z)$ for any $z\in\HH$ and $\gamma\in\Gamma$ of the form $\gamma z= \dfrac{\mathscr az+\mathscr b}{\mathscr cz+\mathscr d}$; this relation has to be valid for some fixed choice of branches of factor $\left(\dfrac{\mathscr cz+\mathscr d}{\mathscr c\bar z+\mathscr d}\right)^\tau$ consistent with the group action; see \cite{Fay} for details. 

In \cite{Ze92}, \cite{D21} the following quantum unique ergodicity result has been proven for integer $\tau_n$'s, but is also true for real ones:

\begin{theorem}
	\label{th:old}
	Let $\Gamma<\Isom^+(\HH)$ be a discrete torsion-free group with a compact fundamental domain $F$, whereas $\tau_1, \tau_2, \dots$ be real numbers.
	
	Suppose that functions $u_n\colon\HH\to\mathbb C$, $n=1,2,\dots$, are such that $u_n$ is a $\tau_n$-form with respect to $\Gamma$, normed as 	$\int_F |u_n|^2\,d\mathcal A_2=2\pi  \mathcal A_2(F)$, and such that $D^{\tau_n} u_n=s_n^2 u_n$ in $\HH$ with $s_n\in\R$ \emph{(}$\mathcal A_2$ denotes hyperbolic area measure on $\HH$\emph{)}.

	If $\tau_n, \dfrac{\tau_n}{s_n}\xrightarrow{n\to\infty}\infty$ 
	then sequence $\{u_n\}_{n=1}^\infty$ is quantum  ergodic. \emph{(}Observables for function  $u_n$ are quantized with Planck constant $\hbar=1/\tau_n$.\emph{)} 
\end{theorem}

In fact, this is a quantization of Furstenberg Theorem on unique ergodicity of horocycle flow over a compact hyperbolic surface (\cite{Furstenberg73}, \cite{Marcus}), up to some calculations on gauge invariance. Thus, Theorems \ref{th:growth} and \ref{th:main} can be applied to functions from Theorem \ref{th:old} and give control on their growth and on the behavior of nodal sets of their complexifications.

\bigskip

Now, let us outline the proof of Theorem \ref{th:growth}. We generally follow Zelditch (\cite{Ze07}). In his paper, he studies similar questions on free-particle quantum ergodic wavefunctions on an arbitrary compact manifold with real-analytic Riemannian metric. Geodesic flow is then instead of horocycle flow. Consequently, instead of horocycle Grauert tube, there is the most usual geodesic Grauert tube defined for any Riemannian manifold. We review Zelditch's argument in the discussion after Proposition \ref{prop:pdo}. 

In our paper, from physicist's viewpoint, we quantize magnetic particle on $\HH$. As well as \cite{Ze07}, our proofs fit into the idea of Boutet de Monvel Theorem. The latter principle is as follows: when we move away from the original real manifold into its complexification, growth of complexified eigenfunction $u^\CC$ is governed by  microlocal distribution of original $u$ in the real part of manifold under consideration. This theorem, in a particular case of Laplacian on a real-analytic manifold,  was stated in \cite{Bou} and has been proved much later in \cite{Ze11}, \cite{Leb13}, \cite{St14}. In this approach, we write $u^{\CC}$ as scalar multiple of $\exp\left(-t\sqrt{-\Delta}\right)u$ and continue analytically Schwartz kernel of $\exp\left(-t\sqrt{-\Delta}\right)$. Studying this operator allows both to provide analytic continuation of $u$ and to set quadratic estimates on $|u^\CC|$.

Unfortunately, our case of horocycle flow and horocycle eigenfunctions  is not covered by the existing results in the spirit of Boutet de Monvel Theorem which concern geodesic setting. In this paper, we first write, using \cite{Fay}, analytic continuation of $u$ to $\mathcal G_{1}$ via an integral operator. For $t\in(0,1)$ understood as horocycle coordinate, we study kernel 
\begin{equation}
\label{eq:kernel_intro}
K_t^\tau(z_1, z_2) = \left(\dfrac{z_1-\bar z_2}{\bar z_1-z_2}\right)^\tau \cdot\exp(-\tau c_t \cosh\dist(z_1,z_2)) ~~~ (z_1, z_2\in\mathbb H)
\end{equation}
with certain $c_t\in\R^+$ (Section \ref{sec:kernel}). Then $A_1u(z_2):=\int_{\HH}K_t^\tau(z_1, z_2) u(z_1) \, d\mathcal A_2(z_1)$ is a scalar multiple of $u$ whenever $u$ is an eigenfunction of $D^\tau$. As it is provided by \cite{Fay}, any kernel of the form
\begin{equation}
\label{eq:general_kernel}
\tilde K(z_1, z_2) = \left(\dfrac{z_1-\bar z_2}{\bar z_1-z_2}\right)^\tau \cdot \left(\mbox{function of }\dist(z_1,z_2)\right)
\end{equation}
has such a property, under reasonable summability conditions.
Term 
\begin{equation}
\label{eq:gauge_factor}
G^\tau(z_1,z_2) := \left(\dfrac{z_1-\bar z_2}{\bar z_1-z_2}\right)^\tau
\end{equation}
is understood as \emph{gauge factor} which also makes these kernels automorphic with respect to isometries of $\HH$. The presence of gauge factor is one of the principal features making our considerations different from that of \cite{Ze07}.

Our kernel is also such that $\HH\ni z\mapsto K_t^\tau(z_1, z)$ can be continued analytically to $\mathcal G_{1}$ (the same concerns mapping $z\mapsto G(z_1,z)$ for gauge factor). This leads to explicit integral formula for $u^\CC$ on $\mathcal G_{1}$. Then, we may put this formula to left-hand side of the limit relation in Theorem \ref{th:growth}. We see that weighted averaging of $|u^\CC|^2$ over $\mathcal G_1$ leads us to a composition of operators in spirit of $A_1^*\mathcal M_a A_1$ acting on functions on $\HH$; 
here $\mathcal M_a$ is multiplication by $a$ acting on functions on $\mathcal G_1$, and now $A_1$ is operator with kernel $K_t^\tau(z_1,z_2)$ continued to $\mathcal G_1$ analytically with respect to $z_2$. 

Fix $t\in(0,1)$. 
Slice $\Sigma_t$ (see (\ref{eq:slice_def})) is homeomorphic to (co)spherical bundle over $\HH$ and thus is naturally endowed with invariant Liouville measure $dS_t$, see Section \ref{sec:coords_flows}. Define diffeomorphism $M_t\colon \{H_1=1/2\}\to \Sigma_t$:
$$
M_t\left(\mbox{covector } \frac{(1+\cos\theta)\,dx+\sin\theta\,dy}y \mbox{ at }x+iy\\
\right) := h_{-it}(x+iy,\theta)
$$
for $x+iy\in\HH$, $\theta\in\R\mmod2\pi$; any point in $\{H_1=1/2\}$ can be parametrized as at the left. 
Operator given by $K_t^\tau(z,P)$ ($z\in\HH$, whereas this time $P\in\Sigma_t$) should be, intuitively and very roughly speaking, understood as \emph{semiclassical \emph{($\hbar=1/\tau$)} Fourier Integral Operator with complex phase and canonical graph}
\begin{multline}
\label{eq:canonical_graph}
\{\left((z,\xi),(M_t(z,\xi), \mbox{some covector at }M_t(z, \xi)\right)\colon (z,\xi)\in\{H_1=1/2\}\}\subset \\ \subset \{H_1=1/2\}\times T^*\Sigma_t\subset T^*\HH\times T^*\Sigma_t.
\end{multline}
To hit the level set $\{H_1=1/2\}$ supporting semiclassical measure of functions $u$, we have to adjust $c_t$, the parameter in kernel $K_t$, see (\ref{eq:kernel_intro}). This "canonical graph"{} lacks dimension, and we repair this by a mollification via $g(\eta)$, see more details in Section \ref{sec:pdo}.

Unfortunately, to author's best knowledge, there is no theory of operators of such a kind.  To calculate a "composition" we apply complex stationary  phase method (\cite{Treves2}, \cite{Hor1}). To this end, we need a global maximum property given by Lemma \ref{lemma:kernel_define_and_max}.


In this manner, in Section \ref{sec:pdo} we construct smooth functions $\mathpzc b(z,\xi)\colon T^*\HH\to (0,+\infty)$, $B(P)\colon \mathcal G_{1}\setminus\HH\to (0,+\infty)$ with the following property. For any $a\in C_0^\infty(\Sigma_t)$, there exists a smooth symbol $\underline{a} \colon T^*\HH\to \R$ such that, first, $\underline{a}$ coincides to $\mathpzc b\cdot (a\circ M_t)$ on $\{H_1=1/2\}$, second, for pseudodifferential operator $A := \Op_{1/\tau} \underline a$, we have 
\begin{equation*}
	\int\limits_{\Sigma_t}dS_t(P)\, a(P) B(P) |u^\CC(P)|^2=\tau^{-3}\cdot \langle Au,u \rangle_{L^2(\HH)}+O(\tau^{-4}) ~~ \mbox{ as } \tau\to+\infty
\end{equation*}
(see Propositions \ref{prop:pdo} and \ref{prop:weakstarconvergence} for more precise statement). Function $B$ is given by an expression depending on $\tau$ and  $s$  but not on $u$.  As we indicated above, auxiliary mollifier $g$ is involved both in $B$ and $\mathpzc b$. Semiclassical measure of sequence $\{u\}$ is concentrated at $\{H_1=1/2\}$. Thus, to calculate the asymptotics of the right-hand side, it is enough to know the symbol of $A$ only on this critical energy level set.
 

To arrive to Theorem \ref{th:growth}, it remains to calculate asymptotics for $B$ as $\tau\to+\infty$ (Section \ref{sec:B_asymp}), here we also eliminate $g$. It requires more applications of stationary phase and Laplace method. Function $B_0$ figuring at the answers in Theorem~\ref{th:growth} and Theorem~\ref{th:main} is obtained in the following manner. Any $P\in\mathcal G_1$ can be written as $P=h_{-it}(z,\theta)$ for some $z\in\HH$, $t\in(0,1)$ and $\theta\in\R\mmod2\pi$. Then
$$ 
B_0(P) = -2\log\left|G(z,P)\right|
$$
with gauge factor $G$ being defined in (\ref{eq:gauge_factor}). We thus may give a brief and qualitative reformulation of Theorem \ref{th:growth}:
$$
\parbox{12cm}{\emph{Growth of a complexified horocycle eigenfunction is given by \\the growth of kernel gauge factor restricted to the canonical graph.}}
$$

\bigskip

\noindent {\bf Acknowledgments.} I am grateful to Steven Zelditch for encouraging to make this paper better. Alas, we could discuss only a very preliminary version of this paper. Since my proofs  generally copy Zelditch's ones, this paper can be considered as a tribute to Steven.

I am also grateful to The Unknown Reviewer for  questions provoked me to improve this paper.

I used \cite{Sage} for most routine symbolic calculations. 

\section{Coordinates and flows}

\label{sec:coords_flows}

In this paper, we denote by $\HH$ the standard upper-halfplane model of Lobachevsky hyperbolic plane. Metric tensor in $\HH$ is given by $(dx^2+dy^2)\cdot y^{-2}$, $x+iy\in \HH$, $x\in\R$, $y>0$.

Point $(X,Y)\in \HH^\CC:=\CC\times\CC$ will be generally denoted by $P$, we write $X(P)$ for $X$ and $Y(P)$ for $Y$. Complex structure in $\HH^\CC$ is that of $\CC^2$. Thus, mappings $P\mapsto X(P)$ and $P\mapsto Y(P)$ are analytic on $\HH^\CC$. We also use $Z(P) := X(P)+iY(P)$ and $\tilde Z(P) := X(P) - i Y(P)$, the analytic continuations of functions $z$ and, respectively, $\bar z$ from $\HH$ to $\HH^\CC$.

Recall that for $z, w\in\HH$ we have, in the hyperbolic metric, $$\dist(z,w) = \arccosh\left(1+\dfrac{|z-w|^2}{2\Im z \, \Im w}\right).$$ Thus, for $z=x+iy\in\HH$ and $P\in\HH^\CC$ ($Y(P)\neq 0$), we may put $$\cosh \dist(z,P) := 1+\dfrac{(x-X(P))^2+(y-Y(P))^2}{2 y Y(P)},$$ and the latter is single-valued function holomorphic with respect to $P$.

Any orientation-preserving isometry of $\HH$ having canonical form $\gamma z= \dfrac{\mathscr az+\mathscr b}{\mathscr cz+\mathscr d}$,  
$\mathscr a,\mathscr b,\mathscr c,\mathscr d\in\R$, $\mathscr a\mathscr d-\mathscr b\mathscr c=1$, $z\in\HH$, can be extended analytically to $\CC\times\CC$, up to possible zeroes in the denominator:
\begin{equation}
	\label{eq:rotation_in_coords}
\CC\times\CC\ni(X,Y)\mapsto\gamma(X,Y) = \left(\dfrac{(\mathscr aX+\mathscr b)(\mathscr cX+\mathscr d)+\mathscr a\mathscr cY^2}{(\mathscr cX+\mathscr d)^2+(\mathscr cY)^2}, \dfrac{Y}{(\mathscr cX+\mathscr d)^2+(\mathscr cY)^2}\right).
\end{equation}
(As it will be seen soon, zeroes at the denominator do not really occur in our considerations if we restrict our interest only to the Grauert tube of horocycle radius~$1$.) We have $Z(\gamma(P)) = \gamma(Z(P))$ and $\tilde Z(\gamma(P)) = \gamma(\tilde Z(P))$. Obviously, such isometries preserve complexified  $\cosh \dist(\cdot, \cdot)$. For the complexification of \emph{gauge factor} $\dfrac{z_1-\bar z_2}{\bar z_1 - z_2}$, the following relation is useful for calculations:
\begin{equation}
	\label{eq:gauge_isometries}
	\dfrac{\gamma z-\gamma \tilde Z}{\gamma \bar z - \gamma Z} = \dfrac{(\mathscr c\bar z+\mathscr d)(\mathscr cZ+\mathscr d)}{(\mathscr cz+\mathscr d)(\mathscr c\tilde Z+\mathscr d)}\cdot \dfrac{z-\tilde Z}{\bar z- Z}.
\end{equation}
We need one more relation. If $\gamma z = \dfrac{\mathscr a_{\gamma}z+\mathscr b_{\gamma}}{\mathscr c_{\gamma}z+\mathscr d_{\gamma}}$, $\gamma^{-1} z = \dfrac{\mathscr a_{\gamma^{-1}}z+\mathscr b_{\gamma^{-1}}}{\mathscr c_{\gamma^{-1}}z+\mathscr d_{\gamma^{-1}}}$ ($z\in\HH$) are isometries written in the canonical form then, for $P\in\HH^{\CC}$, we have
\begin{equation}
\label{eq:gauge_invert}
\dfrac{\mathscr c_\gamma Z(P)+\mathscr d_\gamma}{\mathscr c_{\gamma}\tilde Z(P)+\mathscr d_\gamma}\cdot \dfrac{\mathscr c_{\gamma^{-1}} Z(\gamma P)+\mathscr d_{\gamma^{-1}}}{\mathscr c_{\gamma^{-1}}\tilde Z(\gamma P)+\mathscr d_{\gamma^{-1}}}=1.
\end{equation}
This is consistent with possibility to put $\gamma^{-1}$ instead of $\gamma$ to the definition of $\tau$-form given before Theorem \ref{th:old} and can be verified directly.

Among all the isometries of $\HH$ we widely use the following two types of them. The first is $z\mapsto y_0z+x_0$ ($z\in\HH$) with $x_0\in\R$, $y_0>0$ fixed. Most of our constructions are obviously invariant with respect to them. The second kind is the set of rotations of $\HH$ around $i$ by some angle $\theta\in\R\mmod2\pi$:
\begin{equation*}
	R_{\theta}z:= \dfrac{z\cos(\theta/2)+\sin(\theta/2)}{-z\sin(\theta/2)+\cos(\theta/2)}.
\end{equation*}

A {(}right{)} \emph{horocycle} on Lobachevsky plane $\mathbb H$ is a parametrized curve of constant geodesic curvature $1$ curving to the right and passed with the \emph{unit} speed. An equivalent definition is: 1. the curve $t\mapsto (-t,1)$, $t\in\R$, in $(x,y)$-coordinates in $\HH$ is a right horocycle, 2. any shift of this curve by an isometry of $\HH$ is also a horocycle. 

We widely use \emph{horocycle coordinates} in subsets in $\HH^\CC$. Let $z=x+iy\in\HH$, $\theta\in\R\mmod 2\pi$, $t\in\R$. Let 
\begin{equation}
	\label{eq:vec_parametrize}
	v={y}\cdot\left(\cos\theta\dfrac{\partial}{\partial x} + \sin\theta\dfrac{\partial}{\partial y}\right)\in T_z\HH
\end{equation}
be unit vector based in $z$. There exists a  unique horocycle parametrized as $t\mapsto \phi(t)$, $t\in\R$, with $\phi'(0)=v$, $\phi(0)=z$. Put $h_t(z,\theta):=\phi(t)\in\HH$. Obviously, $\Re h_t(z,\theta)$, $\Im h_t(z,\theta)$ depend analytically on $t$. Therefore, mapping $t\mapsto h_t(z, \theta)$ with $z$, $\theta$ fixed admits an analytic by $t$ continuation for complex $t$ near $\R$.

More precisely, let $t\in\R$. If $\theta=\pi$ then $h_{t}(x+iy,\pi)=x-ty+iy$, $t\in\R$. Thus, 
\begin{equation}
	\label{eq:horizontal_horocycle}
h_{-it}(x+iy,\pi)=(x+ity,y)\in \CC\times\CC
\end{equation}
is its complexification. This "horizontal"{} horocycle is the simplest one, and, if we have some rotation invariance then we often prefer to make calculations in the case $\theta=\pi$, $x+iy=i$.
 
Apply inversion $z\mapsto -1/z$ to the real-time horocycle $h_{t}(x+iy,\pi)=x-ty+iy$. We see that any other non-horizontal right horocycle can be, up to time shift, parametrized as 
\begin{equation}
	\label{eq:horocircle_parameter}
	t\mapsto x_0+y_0\cdot \dfrac{1}{t-i}=x_0+\dfrac{y_0 t}{t^2+1}+i\cdot \dfrac{y_0}{t^2+1}
\end{equation}
with some $x_0\in\R$, $y_0\in\R^+$. We may put complex $t$ with $|\Im t|<1$ to the real and imaginary parts of the latter formula. 
Form  parametrizations (\ref{eq:horizontal_horocycle}) and (\ref{eq:horocircle_parameter}), we derive the following 

\begin{predl}
	For $|\Im t|<1$, mapping $t\mapsto h_{t}(z, \theta)\in\HH^\CC$ can be defined correctly such that $X(h_{t}(z,\theta))$, $Y(h_{t}(z,\theta))$ depend analytically on~$t$ when $z=x+iy\in\HH$ and $\theta \in \R\mmod2\pi$ are fixed.
\end{predl}

%


Further, we have 

\begin{predl}
	\label{predl:Grauert}
	Mapping $\HH\times(0,1)\times(\R\mmod 2\pi)\ni (z,t,\theta)\mapsto h_{-it}(z,\theta)\in \HH^\CC$ is a diffeomorphism onto a set of the form $U\setminus \HH$ with $U\subset\HH^\CC$ being an open neighbourhood of $\HH$. 
\end{predl}

\begin{define}
	\label{def:Grauert_sets}
	For $\tilde t\in (0,1)$, set 
	$$
	\mathcal G_{\tilde  t} := \{h_{-it}(z,\theta)\colon t\in [0,\tilde  t),\, z\in\HH, \, \theta\in\R\mmod 2\pi\}\subset\HH^\CC
	$$
	is called \emph{horocycle Grauert tube} of radius $\tilde t$. \emph{(}Notice that we may take $t=0$ and thus $\HH\subset \mathcal G_{\tilde  t}$ for any $\tilde  t$.\emph{)}
	
	Define also \emph{slice} $\Sigma_{\tilde  t} := \{h_{-i\tilde  t}(z,\theta)\colon z\in\HH, \, \theta\in \R\}$, this is the boundary of $\mathcal G_{\tilde t}$, and, for $\theta\in\R\mmod2\pi$, put $\Sigma_{\tilde  t, \theta} := \{h_{-i\tilde  t}(z,\theta)\colon z\in\HH\}$.
\end{define}

Notice that $\mathcal G_{t_1}\subset \mathcal G_{t_2}$ for $t_1<t_2$ and $\bigcap\limits_{t\mbox{ \scriptsize small}}\mathcal G_{t}=\HH$. Of course, factor $(0,1)\times(\R\mmod2\pi)$ in the domain of mapping in Proposition \ref{predl:Grauert} should be understood as a punctured disk, so that any $\mathcal G_{t}\setminus\HH$, $t\in(0,1)$, is homeomorphic to $\HH\times(\mbox{punctured disk})$. The set of the latter punctures is $\HH$. Thus, we may think about $\mathcal G_t$ as about $\HH\times(\mbox{ball})$, the (co)ball bundle over $\HH$.

We will see soon that horocycle Grauert tubes coincide to geodesic Grauert tubes. Despite this, we call them {horocycle} tubes since we work with horocycle parametrization (as in Proposition \ref{predl:Grauert}) of these sets.

\medskip

\noindent {\bf Proof of Proposition \ref{predl:Grauert}.} For  $x+iy\in\HH$, put 
\begin{gather*}
l^{(1)}_{x+i y}:=\left\{\left(x-\dfrac{y t}{t^2+1}, \dfrac{y}{t^2+1}\right)\colon t\in \CC, \, -1<\Im t<0\right\}\subset \HH^\CC,\\ 
l^{(2)}_{x+i y}:=\{(x+i y t, y)\colon  t\in \R, \, 0< t<1\}\subset\HH^\CC
\end{gather*}
(see (\ref{eq:horizontal_horocycle}), (\ref{eq:horocircle_parameter}) which indeed do parametrize all the horocycles, either non-horizontal or horizontal ones, respectively). 
To prove injectivity from our statement it is enough to show that any of two sets of the form $l^{(1)}_{x+i y}$, $l^{(2)}_{x+i y}$ are disjoint when $x, y$ are varying. We consider the case of two sets of the first kind, the other cases are simpler. Suppose that  $x\in\R$, $y>0$, $t=t_1+i t_2$ ($t_1\in\R$, $t_2\in(0,1)$), $X=X_1+iX_2$, $Y=Y_1+i Y_2$  ($X_1, X_2, Y_1, Y_2\in\R$) and $x-\dfrac{y t}{t^2+1}=X$, $\dfrac{y}{t^2+1}=Y$. Then $X+tY\in \R$, $X_2+t_2 Y_1+t_1 Y_2=0$, also $Y(1+t^2)\in\R$ and  $Y_2(t_1^2-t_2^2+1)+2Y_1t_1t_2=0$. Substituting 
\begin{equation}
\label{eq:t_1_Y_2}
t_1=-(X_2+t_2 Y_1)/Y_2
\end{equation}
to the latter, we find $t_2^2=(X_2^2+Y_2^2)/(Y_1^2+Y_2^2)$ which allows to recover $t_2$ from $X$ and $Y$. (The case $Y_2=0$ is simpler.) Then $t_1, y$ and $x$ are also defined uniquely by $X$ and~$Y$. 

In the remaining cases we also have injectivity. (By the way, we may just rotate $\HH$ to avoid the case of horizontal horocycles, the same concerns the case $Y_2=0$ above.) Now, let us prove that $\mathcal G_1$ contains a neighbourhood of $\HH$ in $\HH^\CC$. Recall that $h_{-it}(i,\pi)=(it, 1)$. Application of complexified rotation by angle $\pi+\theta$ 
around $i$ and also of mapping $z\mapsto x_0+y_0\cdot z$, $z\in\HH$, for fixed $x_0+iy_0\in\HH$ lead to coordinate expressions for $(X,Y)=h_{-it}(x+iy,\theta)$:
\begin{equation}
\begin{gathered}
	\label{eq:cumbersome}
	\Re X = x+y\cdot \dfrac{\left(t^4-2t^2+(t^4-4t^2)\cos\theta\right)\sin\theta}{t^4+(t^4-4t^2)\cos^2\theta+2(t^4-2t^2)\cos\theta+4},\\
	\Im X =y\cdot \dfrac{2\left(t^3+(t^3-2t)\cos\theta\right)}{t^4+(t^4-4t^2)\cos^2\theta+2(t^4-2t^2)\cos\theta+4},\\
	\Re Y =y\cdot \dfrac{2\left(2-(1+\cos\theta)t^2\right)}{t^4+(t^4-4t^2)\cos^2\theta+2(t^4-2t^2)\cos\theta+4},\\
	\Im Y =-y\cdot \dfrac{4t\sin\theta}{t^4+(t^4-4t^2)\cos^2\theta+2(t^4-2t^2)\cos\theta+4}.
\end{gathered}
\end{equation}
From these expressions one can derive the following: if we take  new variables  $x,y$, $v_1=t\cos\theta$, $v_2=t\sin\theta$ such that $(x,y,v_1\frac{\partial}{\partial x}+v_2\frac{\partial}{\partial y})$ runs $T\HH$, then mapping $(x,y,v_1, v_2)\mapsto (X,Y)$ is $C^1$-smooth and  
$$
\left.\dfrac{\partial(\Re X, \Im X, \Re Y, \Im Y)}{\partial(x,y,v_1, v_2)}\right|_{v_1 = v_2=0} = 
\begin{pmatrix}
	1 & 0 & 0 & 0\\
	0 & 0 & -y & 0\\
	0 & 1 & 0 & 0\\
	0 & 0 & 0 & -y
\end{pmatrix}.
$$
This matrix is non-degenerate. It follows that $\mathcal G_1$ indeed contains a neighbourhood of~$\HH$.

It remains to show that Jacobian of 
$$
\det\dfrac{\partial(\Re X, \Im X, \Re Y, \Im Y)}{\partial(x,y,t, \theta)} ~ \mbox{ with }(X,Y)=h_{-it}(x+iy,\theta)
$$ 
is non-zero when $t>0$. Applying isometry of $\HH$ we may assume that 
$\theta=\pi$. Using (\ref{eq:cumbersome}) we find:
\begin{equation}
	\label{eq:Jacobi}
\left.\dfrac{\partial(\Re X, \Im X, \Re Y, \Im Y)}{\partial(x,y,t, \theta)}\right|_{\theta=\pi} = 
\begin{pmatrix}
	1 & 0 & 0 & -\frac{t^2y}2\\
	0 & t & y & 0\\
	0 & 1 & 0 & 0\\
	0 & 0 & 0 & ty
\end{pmatrix}.
\end{equation}
This matrix is, indeed, non-degenerate.
Proof is complete.
$\blacksquare$

\medskip

\noindent {\bf Remark.} From the proof we observe that 
\begin{equation}
\label{eq:cute_formula}
t=\sqrt{\dfrac{(\Im X)^2+(\Im Y)^2}{(\Re Y)^2+(\Im Y)^2}}.
\end{equation}
It follows, in particular, that the expression under square root is invariant with respect to  M\"obius isomorphisms extended to $\mathcal G_1$. Also, (\ref{eq:cute_formula}) or  (\ref{eq:horocircle_parameter}) imply that 
\begin{equation}
	\label{eq:ReY_ImX}
	\Re Y > |\Im X| 
\end{equation}
on $\mathcal G_1$ provided that $t<1$. Proceeding calculations from the proof of Proposition~\ref{predl:Grauert} we see that it is indeed enough for $(X,Y)$ to belong to $\mathcal G_1$: condition (\ref{eq:ReY_ImX}) is also implies $y>0$. Thus, $\mathcal G_1=\{(X,Y)\in\CC^2\colon \Re Y>|\Im X|\}$. From this, one may also conclude the following: if $(X,Y)$ ranges $\mathcal G_1$ then $(Z, \tilde Z)$ ranges $\CC^+\times \CC^-$ (where $\CC^-:=\{z\in\CC\colon \Im z <0\})$. Moreover, action of the group of complexified M\"obius transforms on $\mathcal G_1$ diagonalizes  in $(Z, \tilde Z)$-chart: if $\gamma^\CC\colon \mathcal G_1\to \mathcal G_1$ is complexification of an isometry $\gamma\colon\HH\to\HH$ then $\gamma^\CC(Z, \tilde Z)=(\gamma Z, \gamma \tilde Z)$ with $\gamma$ defined on $\CC^-$ by the same rational fraction as on $\CC^+$.

%
%

If, almost as in the latter proof, we put $v_1=yt\cos\theta$, $v_2=yt\sin\theta$ then, by (\ref{eq:cumbersome}), $X$ and $Y$ do \emph{not} depend $C^\infty$-smoothly on $v_1, v_2$ near $v_1=v_2=0$ due to lots of terms of the kind $t^{\mbox{\scriptsize{odd power}}}$, $t^{\mbox{\scriptsize{even power}}}\cdot \cos\theta$, $t^{\mbox{\scriptsize{even power}}}\cdot \sin\theta$ in (\ref{eq:cumbersome}). But, as we have seen at the proof, parametrization $\left(x+iy,v_1\frac{\partial}{\partial x}+v_2\frac{\partial}{\partial y}\right)\mapsto h_{-it}(x+iy,\theta)$ is $C^1$-smooth up to $t=0$ where it parametrizes $\HH\subset\mathcal G_1$. Also, we may proceed calculations from the proof of Proposition~\ref{predl:Grauert} to ensure that the latter parametrization is proper map from open unit ball bundle over $\HH$ onto $\mathcal G_1$.

%

Notice also that $\{\Im Y=0\}\cap \mathcal G_1=\{\theta\in\{0,\pi\}\}\sqcup \HH$, this is seen from our parametrizations of horocycles. 
\medskip

Hyperbolic plane $\HH$ is endowed with Riemann area $d\mathcal A_2 = \dfrac{dx\, dy}{y^2}$. Tangent spherical bundle $\mathscr S_1\HH$ is endowed with Liouville measure $\tilde\mu_L$: if vectors from this bundle are parametrized as in (\ref{eq:vec_parametrize}) then $d\tilde\mu_L = \dfrac{dx\, dy\, d\theta}{y^2}=d\mathcal A_2 \, d\theta$.

Function $H_1$ defined at Introduction and understood as a Hamiltonian generates bijective identification $\psi_1\colon T\HH\to T^*\HH$ given by 
$$
\psi_1\left(x,y,v_x\frac{\partial}{\partial x}+v_y\frac{\partial}{\partial y}\right) = \left(x,y,\frac{v_x+y}{y^2}dx+\frac{v_y}{y^2}dy\right)\in T^*_{x+iy}\HH
$$
for $x+iy\in \CC^+$, $v_x\frac{\partial}{\partial x}+v_y\frac{\partial}{\partial y}\in T_{x+iy}\HH$
(see also \cite{Takhtajan}). Horocycle Liouville measure $\mu_L$ on the set $\{H_1=1/2\}$ mentioned in the Introduction is given by $\mu_L := (\psi_1)_\sharp \tilde{\mu}_L$, this is the  push-forward of $\tilde\mu_L$ by mapping $\psi_1$.

Push-forward of measure $\tilde\mu_L$ by the mapping
$$
\mathscr S_1\HH \cap T_{x+iy}\HH \ni {y}\cdot\left(\cos\theta\dfrac{\partial}{\partial x} + \sin\theta\dfrac{\partial}{\partial y}\right) \mapsto h_{-it}(x+iy, \theta) \in \Sigma_t
$$
will be denoted by $S_t$, this is the uniform measure on slice $\Sigma_t$.

\medskip

Now, let us compare horocycle Grauert tubes to the usual \emph{geodesic} Grauert tubes (see \cite{GS91}, \cite{GS92}, \cite{LS91}): as sets (endowed with complex structure!) they coincide up to change of imaginary time. Denote by $h^{\geod}_r(z,\theta)$ the unit-speed geodesic line starting 
in $z\in\HH$ with the slope $\theta\in\R$ to the real axis and evaluated at the time $r\in\R$. We may complexify it by $r$ and consider $h^{\geod}_{ir}(z,\theta)\in\HH^\CC$ for $r$ small. The above-mentioned geodesic Grauert tube is $\mathcal G_r^{\geod}=\{h^{\geod}_{ir_1}(z,\theta)\colon z\in\HH, \,\theta\in\R, \, r_1\in[0,r)\}$ where $r$ again is 
not too large. Introduce also slices $\Sigma_r^{\geod}:=\{h^{\geod}_{ir}(z,\theta)\colon z\in\HH, \,\theta\in\R\}$.

We have $h^{\geod}_r(i, \pi/2)=ie^r\in \HH$. Thus, $h^{\geod}_{ir}(i, \pi/2)=(0, e^{ir})\in\CC^2$. Applying 
an isometry  (\ref{eq:rotation_in_coords}) 
we see that zero in the denominator arises when $e^{2ir}=-1$. Thus, $r$ may range $(0,\pi/2)$. 

Observe that the group of complexified isometries of $\HH$ preserves  any of the sets $\Sigma_t$ or $\Sigma_r^{\geod}$  for  fixed $t\in(0,1)$ or, respectively, $r\in(0,\pi/2)$, and acts transitively on any of them. It follows that any two of such sets either coincide or do not intersect. 
If $r\in(0,\pi/2)$ then $h^{\geod}_r(i,\pi/2)\in \{\Re Y>|\Im X|\}$ which is $\mathcal G_1$; by invariance of $\mathcal G_1$ with respect to the group of complexified isometries of $\HH$, we conclude that $\mathcal G^{\geod}_{\pi/2}\subset \mathcal G_1$. To show that 
$\mathcal G_1\subset \mathcal G^{\geod}_{\pi/2}$
pick any $t\in(0,1)$. We have  $h^{\geod}_{\arcsin t}(i,\pi/2)\in \Sigma_t$ and therefore $\Sigma^{\geod}_{\arcsin t}=\Sigma_t$. It follows that $\mathcal G_1= \mathcal G^{\geod}_{\pi/2}$.

Slice-wise recalculation from horocycle to geodesic coordinates $(x,y,\theta)$ in $\Sigma_t=\Sigma^{\geod}_{\arcsin t}$ is non-degenerate since it can be reduced to group multiplication in $\Isom^+(\HH)$, the group of orientation-preserving isometries of $\HH$, which is identified to $\mathscr S_1\HH$, the spherical bundle over $\HH$. This implies non-degeneracy of geodesic parametrization for $r\in(0,\pi/2)$. We arrive to the following conclusion: 


\begin{predl}
\label{predl:Grauert_hor_geod} In half-plane hyperbolic plane model, maximal radius of geodesic Grauert tube is $\pi/2$. We have
$\mathcal G_1=\mathcal G_{\pi/2}^{\geod}$, and, for $t\in(0,1)$, we have $\Sigma_t= \Sigma^{\geod}_{\arcsin t}$ as sets with complex \emph{(}respectively, CR-\emph{)} structure.
\end{predl}

\begin{sled}
\label{sled:involution} Any of sets $\Sigma_t$, $t\in(0,1)$, is invariant with respect to the  involution $\imath\colon \HH^\CC\to\HH^\CC$ defined as $\imath(X,Y):= (\bar X, \bar Y)$.
\end{sled}

This is because it is true for $\Sigma^{\geod}_{\arcsin t}$. Also, this follows from the Remark after Proposition \ref{predl:Grauert}.

In the usual geodesic Grauert tube endowed with imaginary geodesic parametrization, there is an intriguing plexus of structures leading to K\"ahler geometry. For $\mathcal{G}_{\pi/2}^{\geod}\ni P = h_{ir}^{\geod}(z,\theta)$ ($z\in\HH$, $r,\theta\in\R$) put $f(P):=r$. Then:
\begin{itemize}
	\item $f^2$ is plurisubharmonic in $\mathcal G_{\pi/2}^{\geod}$ with respect to the complex structure in $\HH^\CC=\CC\times\CC$;
	
	\item $\displaystyle\det\left(\dfrac{\partial^2}{\partial V \partial  \overline{W}}f\right)_{V,W\in\{X,Y\}}=0$ out of $\HH$ (complex Monge--Amp\`ere equation);
	
	\item form $-i\partial\bar{\partial} f^2$ 
	is  a symplectic form $\omega_{\mathrm{symp}}$ on $\mathcal G_{\pi/2}^{\geod}$. Let us identify $$(z,r,\theta)\in\HH\times[0,\pi/2)\times (\R\mmod 2\pi)$$ with covector $y^{-1}r\cos\theta\,dx+y^{-1}r\sin\theta\,dy\in T^*_z\HH$. Then, under parametrization $h^{\geod}_{ir}(z,\theta)$ of $\mathcal G^{\geod}_{\pi/2}$, form $\omega_{\mathrm{symp}}$ is taken to $d\xi_1\wedge dx+d\xi_2\wedge dy$, the standard symplectic form at $T^*\HH$.
	
	\item Form $\omega_{\mathrm{symp}}(\cdot, J\cdot)$ is positively defined where $J$ is the complex structure in $\CC\times\CC\supset\mathcal G_{\pi/2}^{\geod}$. 
	
	Thus, $\mathcal G_{\pi/2}^{\geod}$ becomes a K\"ahler manifold.
\end{itemize}
Similar assertions hold if we start with an arbitrary real-analytic Riemannian manifold. We proceed the discussion on waves in geodesic Grauert tubes in Section \ref{sec:pdo}.

Alas, we think that most of this geometric harmony crushes in the horocycle setting.

To this end, return to the horocycle coordinates $(z,t,\theta)$. Let us say that  a $2$-form on $\mathcal G_1\setminus\HH$ is \emph{fiber-orthogonal} if its coefficient before $d\theta\wedge dt$ in $(z,t,\theta)$-chart is zero. This means that our form vanishes at any $2$-vector $F\cdot \frac{\partial}{\partial t}\wedge\frac{\partial}{\partial\theta}$ in the same coordinates.

\begin{predl}
\label{predl:no_fiber}
Let $E=E(t,\theta)\colon \mathcal G_1\setminus\HH\to \CC$ be a function not depending on $z\mbox{-coordinate}$. If $\partial\bar\partial E$ is fiber-orthogonal then $E$ is of the form $f(\theta)+\dfrac{\const}{t}$. 

For such $E$, $\partial\bar\partial E$ is either zero or, in $(X,Y)$-chart, it has singularity near $t=0$; more precisely, in this case,  coefficients of $\partial\bar{\partial}E$ near $\HH$ in the latter chart are not bounded.
\end{predl}

We parametrize $\mathscr B_1(\HH)$, the radius $1$ open ball bundle over $\HH$, by $(z,t,\theta)$-chart sending $(t,\theta)$ to vector $yt\cos\theta\frac{\partial}{\partial x}+yt\sin\theta\frac{\partial}{\partial y}\in T_z\HH$. Any smooth Hamiltonian identification of $T^*\HH$ and $T\HH$ preserves fibers  and thus takes form $d\xi_1\wedge dx+d\xi_2\wedge dy$ to a fiber-orthogonal form in $(z,t,\theta)$-chart. The latter form cannot have a singularity near $t=0$ since parametrization 
$$
\mathscr B_1(\HH)\ni \mbox{ vector }
yt\cos\theta\frac{\partial}{\partial x}+yt\sin\theta\frac{\partial}{\partial y}\mbox{ at point }z \mapsto h_{-it}(z,\theta)\in\mathcal G_1
$$
is $C^1$-smooth up to $\{t=0\}$, as we have seen in the proof of Proposition \ref{predl:Grauert}.  We conclude that functions $E\colon \mathcal G_1\to\CC$ 
depending only on $t$ and $\theta$ cannot have $\partial\bar\partial E$ equal to the standard symplectic $2$-form on $T^*\HH$ transferred smoothly fiberwise to $T\HH$ and, further, to $\mathcal G_1$ via horocycle parametrization.

\medskip

\noindent {\bf Proof of Proposition \ref{predl:no_fiber}.} To avoid second-order complex differentiation in chart $(z,t,\theta)$, we use the identity $\partial\bar\partial=-\frac i2dd^c$, $d^c=i(\bar{\partial}-\partial)$. We calculate forms at the basepoint $z=i$ and with $\theta=\pi$, and then rotate them. Using matrix (\ref{eq:Jacobi}) and its inverse we find, at $\theta=\pi$, that
$$
-{d^cE} =  E'_t t^2d\theta +2E'_t\cdot\frac{dx}{y}+\frac{2E'_\theta}{t}\cdot\frac{dy}{y}.
$$
For $\theta\neq \pi$, we make conformal rotation around basepoint; it follows that ${d^cE}$ is always of the form $-{E'_t t^2} d\theta+A_1dx+A_2dy$. If its differential is fiber-orthogonal then $E'_tt^2=\const$ which implies that $E$ is $C/t+f(\theta)$, $C\in\CC$, $f\colon \R\mmod2\pi\to\CC$. 

Now we proceed routine calculations for this $E$ and arrive to 
$$
\left.-dd^{c}E\left[\frac{\partial}{\partial x}\wedge\frac{\partial}{\partial y}\right]\right|_{z=i}=
-C\cdot\left(\cos\theta+1\right)+\frac{2C\cos\theta}{t^2}+\frac{2f'(\theta)\cdot \sin\theta}{t},
$$
here square brackets denote application of a $2$-covector to a $2$-vector. But $\frac{\partial}{\partial x}\wedge\frac{\partial}{\partial y}$ is bounded in $(X,Y)$ chart for $t$ close enough to zero. This leads to the desired.
$\blacksquare$

\medskip

Now we pass to structures at odd-dimensional $\Sigma_t$. 
At $\theta=\pi$, complex structure in $\HH^\CC$ can be expressed in coordinate frame $\left(\dfrac{\partial}{\partial x},\dfrac{\partial}{\partial y},\dfrac{\partial}{\partial t},\dfrac{\partial}{\partial\theta}\right)$ as 
$$
\begin{pmatrix}
	0 & -\frac{t}{2} & -y & 0 \\
	0 & 0 & 0 & -ty \\
	\frac{1}{y} & 0 & 0 & \frac{t^2}{2}  \\
	0 & \frac{1}{ty} & 0 & 0
\end{pmatrix}.
$$
In terms of \emph{real} tangent bundle $T\Sigma_t$, CR-structure on $\Sigma_t$ is given by subbundle spanned by $\left\{\dfrac{\partial}{\partial y},J\dfrac{\partial}{\partial y}\right\}$, $J\dfrac{\partial}{\partial y}=-\dfrac{t}{2}\cdot\dfrac{\partial}{\partial x}+\dfrac{1}{ty}\cdot\dfrac{\partial}{\partial\theta}$. By Proposition \ref{predl:Grauert_hor_geod}, this $2$-dimensional distribution does not depend on whether we choose horocycle or geodesic parametrization for $\Sigma_t$.

Next, we want to transfer contact structure from $\{H_1=1/2\}$ 
to $\Sigma_t$. We identify $(x+iy,\theta)$ to covector $y^{-1}(\cos\theta+1)dx+ y^{-1}\sin\theta \, dy$ at $x+iy$. Tautological $1$-form $\zeta=\xi_1 dx+\xi_2 dy$ degenerates at $\theta=\pi$ and thus 
cannot be understood as a  contact form, unlike the geodesic setting. Also, $\zeta$ is not rotation-invariant whereas CR-distribution is. Intermediate steps from the proofs from \cite{LS91}, \cite{GS91} concerning $1$-forms thus fail in our context.

It remains to check whether symplectic $2$-form $\omega_{\mathrm{symp}}=d\zeta$ transferred from $\{H_1=1/2\}$ to $\Sigma_t$ is compatible to $J$, that is, form $\omega_{\mathrm{symp}}(\cdot, J\cdot)$ does not change sign. At $\theta=\pi$, $\omega_{\mathrm{symp}}\left(\dfrac{\partial}{\partial y},J\dfrac{\partial}{\partial y}\right)=\dfrac{1}{ty^2}$.

Complex structure $J$ on $\mathcal G_1$ is invariant with respect to rotations as well as its restriction to CR-distribution $\Sigma_t$ since 
all rotations of $\mathcal G_1$ are complex-analytic. Form $\dfrac{\cos\theta\,dx+\sin\theta\,dy}{y}$ also is preserved under rotations. This is because it reduces to  tautological form $\zeta$ transferred to unit spherical bundle $\mathscr S_1\HH$ via \emph{Riemannian} identification of $T\HH$ and $T^*\HH$, and Riemannian structure is preserved under rotations as well as $\zeta$. Finally, $d\left(\dfrac{dx}{y}\right)$ is volume form on $\HH$ and is also invariant with respect to rotations. Thus, $\omega_{\mathrm{symp}}$ is also invariant, and quadratic form $\omega_{\mathrm{symp}}(\cdot, J\cdot)$ is always positive at CR-distribution on $\Sigma_t$.

Besides invariance as just above, we suspect that Theorem on uniqueness of adapted complex structure (\cite[Proposition 5.1]{LS91}) still survives in our setting since its proof seems to be not sensible to the choice of a parametrization.

\section{Construction of an automorphic kernel}

\label{sec:kernel}

In this Section we construct analytic continuation of $u$ to the horocycle Grauert tube $\mathcal G_1$ via an integral operator with kernel $K_t^\tau(\cdot,\cdot)$. For studying the growth of $u^\CC$ at $\Sigma_t$ ($t\in(0,1)$) we need a global maximum property of this kernel. 

\begin{lemma}[on global maximum of absolute value]
	\label{lemma:kernel_define_and_max}
	For $t\in(0,1)$, 
	put $c_t := \dfrac{4}{4t-t^3}$ and 
	\begin{equation}
		\label{eq:half_kernel}
		K_t(z,P) := \left(\dfrac{z-\tilde Z(P)}{\bar z - Z(P)}\right) e^{-c_t \cdot \cosh \dist(z,P)}, ~~ z\in \HH, \, P\in \mathcal G_{1}.
	\end{equation}
\begin{enumerate}
	\item Function $\Phi_t(z,P) = \log K_t(z,P)$ is single-valued when $z\in \HH$ and $P\in \mathcal G_{1}$.
	
	\item For $z\in\HH$ and $\theta\in\R\mmod 2\pi$ fixed,  
	$
	\max\limits_{P\in\Sigma_{t, \theta}}|K_t(z,P)|
	$
	is attained at $P=h_{-it}(z, \theta)$.
	
	\item Hesse matrix
	\begin{equation*}
		\begin{pmatrix}
			\partial_{xx} \Re\Phi_t(z_0, h_{-it}(x+iy, \theta)) & \partial_{xy} \Re\Phi_t(z_0, h_{-it}(x+iy, \theta)) \\
			\partial_{xy} \Re\Phi_t(z_0, h_{-it}(x+iy, \theta)) & \partial_{yy} \Re\Phi_t(z_0, h_{-it}(x+iy, \theta))
		\end{pmatrix}
	\end{equation*}
	is non-degenerate when $x+iy=z_0$.
\end{enumerate}
\end{lemma}	

\noindent {\bf Proof.} First claim follows from (\ref{eq:ReY_ImX}): it implies that ${z_0-\tilde Z(P)}$, ${\bar z_0 - Z(P)}$ are non-zero. 

Now we prove the second assertion. We claim that it is enough to check that 
\begin{equation}
	\label{eq:K_pi_max}
	\argmax{z\in\HH}\log |K_t(z, h_{-it}(i,0))|=i.
\end{equation}
Let us do this reduction. 
First, we claim that (\ref{eq:K_pi_max}) implies that 
\begin{equation}
	\label{eq:K_theta_max}
	\argmax{z\in\HH}\log |K_t(z, h_{-it}(i,\theta))|=i
\end{equation}
for any $\theta\in\R\mmod2\pi$. 
Indeed, for such $\theta$, write $R_{\theta}$ defined at Section \ref{sec:coords_flows} in the canonical form as $\HH\ni z\mapsto R_{\theta}z = \dfrac{\mathscr az+\mathscr b}{\mathscr cz+\mathscr d}$, $\mathscr a,\mathscr b,\mathscr c,\mathscr d\in\R$, $\mathscr a\mathscr d-\mathscr b\mathscr c=1$. If (\ref{eq:K_pi_max}) is already checked then, by (\ref{eq:gauge_isometries}), we have, for any $z\in\HH$, that
\begin{multline*}
|K_t(z, h_{-it}(i,\theta))|=|K_t(R_\theta R_{-\theta}z, R_\theta R_{-\theta}h_{-it}(i,\theta))|=\\=\left|\dfrac{cR_{-\theta}\bar z+\mathscr d}{\mathscr cR_{-\theta}z+\mathscr d}\right|\cdot \left|\dfrac{\mathscr c Z(R_{-\theta}h_{-it}(i,\theta))+\mathscr d}{\mathscr c\tilde Z(R_{-\theta}h_{-it}(i,\theta))+\mathscr d}\right|\cdot |K_t(R_{-\theta}z, h_{-it}(i,0))|\le\\ \le  \left|\dfrac{\mathscr c Z(h_{-it}(i,0))+\mathscr d}{\mathscr c\tilde Z(h_{-it}(i,0))+\mathscr d}\right|\cdot |K_t(i, h_{-it}(i,0))|=\\=|K_t(R_{\theta}i, R_{\theta}h_{-it}(i,0))|=|K_t(i, h_{-it}(i,\theta))|,
\end{multline*}
the desired. 

Now, having (\ref{eq:K_theta_max}), take any $z_0=x_0+iy_0\in\HH$, put $\gamma z:= x_0+y_0\cdot z$ for $z\in\HH$. Notice that 
\begin{multline*}
	|K_t(i, h_{-it}(z,\theta))| = |K_t(\gamma^{-1}i, \gamma^{-1}h_{-it}(z,\theta))|=|K_t(\gamma^{-1}i, h_{-it}(\gamma^{-1}z,\theta))|=\\=|K_t(\gamma^{-1}i, h_{-it}(i,\theta))| \le |K_t(i, h_{-it}(i,\theta))|
\end{multline*} 
by (\ref{eq:K_theta_max}). A similar application of an isometry like the latter one allows to replace $z=i$ in the second assertion of our Lemma by any other point in $\HH$.


Now,  denote 
$P_0:=h_{-it}(i,0)=\left(-\dfrac{it}{1-t^2}, \dfrac{1}{1-t^2}\right)$. 
To check (\ref{eq:K_pi_max}) it is enough to show that $|K_t(x+iy,P_0)|$ with $y$ fixed decreases by $x$ for $x\ge 0$ and increases by $x$ for $x\le 0$, and also that  $|K_t(iy,P_0)|$ attains its maximum over $y\in\R^+$ at $iy=i$. 

We proceed via routine calculation. We have
\begin{multline*}
\frac{\partial\Re\Phi_t(x+iy,P_0)}{\partial x}=\\=
-4x\cdot\frac{A}{\left((x^2+y^2)\cdot(1-t)^2+2y(1-t)+1)
	\cdot((x^2+y^2)\cdot(1+t)^2+2y(1+t)+1)
	\right)(4-t^2)ty}
\end{multline*}
with numerator
\begin{multline*}
A = (1-t^2)^3x^4+(1-t^2)^3y^4+4(1-t^2)^2y^3+2(1-t^4)x^2+(1-t^2)(6+2t^2-t^4)y^2
+\\+
2(1-t^2)^3x^2y^2+
(4-t^4+4(1-t^2)^2x^2)y+1-t^2,
\end{multline*}
which is positive.  Now, compute 
\begin{equation*}
\left.\frac{\partial\Re\Phi_t(x+iy,P_0)}{\partial y}\right|_{x=0}=
2(1-y)\cdot\frac{y^3(1-t^2)^2+y^2(1-t^2)(3-t^2)+3y+1}{ty^2(2-t)(2+t)(1+y+ty)(1+y-ty)}.
\end{equation*}
Numerator at the right-hand side is positive and we arrive to the desired monotonicity.

To prove the third claim, we note that Hesse matrix mentioned therein does not depend on $x+iy$ and, at $x+iy=i$, is the same as 
\begin{equation*}
	\begin{pmatrix}
		\partial_{xx} \Re\Phi_t(x+iy, h_{-it}(i, \theta)) & \partial_{xy} \Re\Phi_t(x+iy, h_{-it}(i, \theta)) \\
		\partial_{xy} \Re\Phi_t(x+iy, h_{-it}(i, \theta)) & \partial_{yy} \Re\Phi_t(x+iy, h_{-it}(i, \theta))
	\end{pmatrix}
\end{equation*}
at $x+iy=i$. By using (\ref{eq:gauge_isometries}) again we see that this matrix does not depend on $\theta$; we thus may take $\theta=0$ as above. Then, by a computation, the latter matrix is
$$
	\begin{pmatrix}
	\dfrac{4 \, {\left(3 \, t^{2} - 4\right)}}{t^{5} - 8 \, t^{3} + 16 \, t} & 0\\
	0 & 	-\dfrac{4 \, {\left(t^{4} - 3 \, t^{2} + 4\right)}}{t^{5} - 8 \, t^{3} + 16 \, t}
	\end{pmatrix}
$$
and is obviously negatively defined as $0<t<1$.
 $\blacksquare$
\medskip

\medskip

Recall that we assume $\tau>0$.

\begin{lemma}
	\label{lemma:SelbergTransform}
	For $u\colon \HH\to \CC$ put $v(P) := \int_{\HH} u(z)K^\tau_t(z,P)\,d\mathcal A_2(z)$ \emph{(}$P\in\mathcal G_{1}$\emph{)}. 
	\begin{enumerate}
		\item If $u\in L^\infty(\HH)$ then the integral above converges absolutely together with any of its derivatives with respect to coordinates of $P$, also uniformly when $P$ ranges a compact set in  Grauert tube $\mathcal G_{1}$ of horocycle radius $1$, and $t$ ranges a compact set in $(0,1)$. 
		
		\item If $u\colon\HH\to\CC$ is such that $D^\tau u =\left(-\Delta_\HH+2i\tau y\dfrac{\partial}{\partial x}\right)u= s^2 u$ then $v(z) = \mathcal S(t,\tau,s) u(z)$ for some $\mathcal S$ not depending neither on $z\in\HH$ nor on $u$.
		\item Function $v(P)$ is analytic for $P\in\mathcal G_{1}$.
	\end{enumerate}

Thus, $v$ is, up to a constant factor, an analytic continuation of $u$ to horocycle Grauert tube $\mathcal G_1$.
\end{lemma}

Proof is rather technical and is given at the Appendix (Subsection \ref{subsec:Sec3proofs}).
%

\medskip

\noindent {\bf Remark.} In Theorem \ref{th:growth} we require, in particular, that $\sup\limits_{z\in \HH} \|u_n\|_{L^1(\mathcal B_\HH(z,1))}<+\infty$ for any given $n$. This implies that each $u_n$ belongs to $L^\infty(\HH)$ (not necessarily uniformly by $n\in \mathbb N$). This can be seen by appropriate averaging the relation from  \cite[Theorem 1.2]{Fay}. Lemma \ref{lemma:SelbergTransform} is therefore applicable to functions from Theorem \ref{th:growth}.

\medskip

\noindent {\bf Remark.} Grauert tube $\mathcal G_1$ is the maximal set to which any eigenfunction can be continued analytically (see also discussion in \cite{Ze11}). 

Indeed, using automorphic change of variables and  (\ref{eq:rotation_in_coords}),  we see that 
$$\HH\ni x+iy=z\mapsto 
\left(\dfrac{\mathscr c\bar z+\mathscr d}{\mathscr c z+\mathscr d}\right)^\tau\cdot\left(\dfrac{y}{(\mathscr cx+\mathscr d)^2+(\mathscr cy)^2}\right)^\alpha$$
is an eigenfunction  of $D^\tau$ for any $\mathscr c, \mathscr d\in\R$, $\alpha\in\mathbb C$ (see \cite[\S1]{Fay}). 
On $\partial\mathcal G_1$, either $Z=X+iY$ or $\tilde Z=X-iY$ is real (see Remark after proof of Proposition \ref{predl:Grauert}). 
It follows that we may pick $\mathscr c$, $\mathscr d$ such that the second factor has singularity at $(X,Y)$; if the first factor also has zero or singularity at this point then it does not cancel the  singularity of the second factor under almost all choices of $\alpha$. Thus, for a given point $P\in\partial\mathcal G_1$, there is a function as above which cannot be analytically continued through~$P$. 
\medskip

Now we perturb $t$ replacing it in $c_t$ but reserving in $\Sigma_{t, \theta}$ (which is, recall, $\{h_{-it}(z,\theta)\colon z\in\HH\}$ with $t,\theta$ fixed). The reason to do so is need to mollify operator with kernel $K_t^\tau$, this is discussed in Introduction and at Section \ref{sec:pdo}. 

\begin{lemma}
\label{lemma:max_point_exists}
	Let $t\in(0,1)$, $\eta$ be close enough to $t$, $\theta\in\R\mmod2\pi$ and $z\in \HH$.
		Function 
	\begin{equation*}
		P\mapsto |K_\eta(z, P)|
	\end{equation*}
	has a unique maximum point at $\Sigma_{t, \theta}$. We denote this maximizer by $Q(z, t, \eta,\theta)$ and also put 
	\begin{equation*}
	\varphi(t,\eta,\theta)	:= \log\max\limits_{P\in\Sigma_{t, \theta}}|K_\eta(z,P)|=\log |K_\eta(z,Q(z,t,\eta,\theta))|.
	\end{equation*}
	
\end{lemma}

\noindent {\bf Proof.} Follows from the third assertion of Lemma \ref{lemma:kernel_define_and_max} and Implicit Function Theorem. 
Also, there is no problem at infinity during this perturbation as can be seen from the proof of the first assertion of Lemma \ref{lemma:SelbergTransform}.
$\blacksquare$ 

\medskip

We need to hit energy level $\{H_1=1/2\}$ in the "canonical graph"{} (\ref{eq:canonical_graph}) since microlocal mass of $u$ (that is, measure as in the relation from Definition \ref{def:QE}) is concentrated near this set. This is provided by the following Lemma~\ref{lemma:diffeo}. Energy level $\{H_{-1}=1/2\}$ obtained in this Lemma will be finally replaced by $\{H_{1}=1/2\}$ by a certain flipping in quadratic form in Proposition~\ref{prop:pdo}.  In the proof of Lemma~\ref{lemma:diffeo}, we, in particular, apply rotations of $\HH^\CC$ around basepoint $z_0$ from the statement. This also allows to calculate the dependence of $\varphi(t,\eta,\theta)$ on $\theta$.


Recall that $\Phi_{\eta}=\log K_\eta$.

\begin{lemma}
\label{lemma:diffeo}
	Let $t\in(0,1)$ be fixed.
	\begin{enumerate}
		\item  For $z_0\in\HH$, mapping $$\mathcal T_{z_0,t}(\eta, \theta):= \Im d_z \left.\Phi_\eta(z, Q(z_0, t, \eta,\theta))\right|_{z=z_0}$$ is a diffeomorphism of $(\mbox{some neighbourhood of }t)\times(\R\mmod 2\pi)$ onto a neighborhood of circle $T^*_{z_0}\HH\cap \{H_{-1}=1/2\}$. 

	\item For any  $\eta$ close enough to $t$, $\theta\in\R\mmod2\pi$, we have, for $B_0$ defined at Introduction,
	$$
	\varphi(t,\eta,\theta) = \varphi(t,\eta,\pi)-\dfrac12\cdot\log\left({\dfrac{(1+\cos\theta)\cdot(t^2-2t)+2}{(1+\cos\theta)\cdot(t^2+2t)+2}}\right) = \varphi(t,\eta,\pi)-\dfrac{B_0}2.
	$$
	\end{enumerate}
\end{lemma}

Clearly, $\mathcal T_{z,t}$ is degree $-1$ homogeneous with respect to $\Im z$. If $\mathcal T_{z,t} (\eta, \theta)=(\xi_1, \xi_2)\in T^*_{z}\HH$ then we write $\theta =: \Theta_{z,t}(\xi_1, \xi_2)$ and $\eta =: \Eta_{z,t}(\xi_1, \xi_2)$. The latter mappings are defined near $\{H_{-1}=1/2\}$. On this level set, if  $z\in\HH$, $\xi_1\,dx+\xi_2\,dy\in T^*_z\HH$, $H_{-1}(z,\xi_1,\xi_2)=1/2$ with $\xi_1=(-1-\cos\theta)/\Im z$, $\xi_2=-\sin\theta/\Im z$, $\theta\in\R$, then  
\begin{equation}
\label{eq:EtaThetaCritical}
\begin{gathered}
	\Theta_{z,t}(\xi_1,\xi_2)=\theta, \\
	\Eta_{z,t}(\xi_1, \xi_2)=t. 
\end{gathered}
\end{equation}
This is seen from the calculations from the proof below.

What concerning the second assertion, the logarithm at the right-hand side will finally lead us to the asymptotics in Theorem \ref{th:growth} and answer in Theorem \ref{th:main}.

\medskip

\noindent {\bf Proof of Lemma \ref{lemma:diffeo}.} Let's start with the first assertion. To begin, we consider the case $\theta=\pi$. 
The point $\argmax{z=x+iy\in\HH} |K_\eta(z,h_{-it}(i,\pi))|$ should be necessarily $iy(t,\eta)$ with some $y(t,\eta)>0$ (the kernel is even with respect to $x$). By routine differentiation, we see that $\left.d_z\right|_{z=i}\Phi_{\eta}(x+iy, Q(i,t,\eta,\pi))=i\cdot f(t,\eta)\,dy$ with $f$ smooth, $f(t,t)=0$ and 
\begin{equation}
		\label{eq:non_degenerate_i_pi}
\left.\dfrac{\partial f}{\partial \eta}\right|_{\eta=t}\neq 0.
\end{equation}

%

To make $\theta\neq \pi$, we apply rotation 
around $i$ by angle $\pi+\theta$. 
Put $R := R_{\pi-\theta}$. Similarly to the second assertion of Lemma \ref{lemma:kernel_define_and_max}, we conclude that 

%

\begin{equation*}
	\label{eq:Q_rorate}
	R^{-1}h_{-it}(i,\pi) = Q(R^{-1}(iy(t, \eta)),t,\eta,\theta)=:Q_\theta.
\end{equation*}
By (\ref{eq:gauge_isometries}),
\begin{multline}
\label{eq:rotate_kernel}
\Phi_\eta(z, Q_\theta) =\\= \Phi_\eta(Rz, RQ_\theta)-\log\left(\dfrac{\bar z\cos(\theta/2)-\sin(\theta/2)}{z\cos(\theta/2)-\sin(\theta/2)}\right)-\log\left(\dfrac{Z(Q_\theta)\cos(\theta/2)-\sin(\theta/2)}{\tilde Z(Q_\theta)\cos(\theta/2)-\sin(\theta/2)}\right)
=\\=
\Phi_\eta(Rz, h_{-it}(i,\pi))+
\log\left(\dfrac{z (\cos\theta+1)-\sin\theta}{\bar z (\cos\theta+1)-\sin\theta}\right)
-\log\left(\dfrac{Z(Q_\theta)\cos(\theta/2)-\sin(\theta/2)}{\tilde Z(Q_\theta)\cos(\theta/2)-\sin(\theta/2)}\right).
\end{multline}
Since $d_z \left.\Phi_\eta(z,h_{-it}(i,\pi))\right|_{z=i}=0$, we have 
\begin{equation}
	\label{eq:dPhi_0_t}
d_z\left.\Phi_t(z,Q_\theta)\right|_{z=i} =  d_z\left.\log\left(\dfrac{z (\cos\theta+1)-\sin\theta}{\bar z (\cos\theta+1)-\sin\theta}\right)\right|_{z=i}=-(1+\cos\theta)\cdot i\, dx-\sin\theta \cdot i\, dy.
\end{equation}

Next, suppose that $F(z)$ is some function and that $\gamma z = \dfrac{\mathscr az+\mathscr b}{\mathscr cz+\mathscr d}$. Suppose that  $\left.d_z F\right|_{z=\gamma z_0} = i\alpha\, dx+i\beta\, dy$ for some $\alpha,\beta\in\mathbb C$. Put  $F_1(z) := F(\gamma z) +\log\left(\dfrac{\mathscr cz+\mathscr d}{\mathscr c\bar z+\mathscr d}\right)$. If $dF_1|_{z_0} = i\alpha_1 \,dx + i\beta_1\, dy$ then
\begin{equation}
	\label{eq:Hamiltonian_invariance}
(\alpha\Im(\gamma z) +1)^2+(\beta\Im(\gamma z ))^2=(\alpha_1\Im z +1)^2+(\beta_1\Im z)^2
\end{equation}
(\emph{energy conservation under automorphic change of variables}).  This  can be checked by a direct calculation.

By (\ref{eq:dPhi_0_t}) we see that
$$
\left.\dfrac{\partial}{\partial\theta}\Im d_z\Phi_\eta(z, Q(z_0,t,\eta,\theta))\right|_{z=z_0, \, \eta=t}
$$
is tangent to $T^*_{z_0}\HH\cap \{H_{-1}=1/2\}$. Further, (\ref{eq:Hamiltonian_invariance}) and (\ref{eq:non_degenerate_i_pi}) imply that 
$$
\left.\dfrac{\partial}{\partial\eta}\Im d_z\Phi_\eta(z, Q(z_0,t,\eta,\theta))\right|_{z=z_0, \, \eta=t}
$$
is transverse to $T^*_{z_0}\HH\cap \{H_{-1}=1/2\}$. Thus, Jacobian of
$$
(\eta,\theta)\mapsto d_z\Im\Phi_\eta(z, Q(z_0,t,\eta,\theta))|_{z=z_0}
$$
is non-zero near ${\eta=t}$. Together with (\ref{eq:dPhi_0_t}), this concludes the proof of the first assertion of our Lemma.

For the second assertion, notice that $Z(h_{-it}(i,\pi))=i(t+1)$, $\tilde Z(h_{-it}(i,\pi))=i(t-1)$. Now apply (\ref{eq:gauge_invert}) with $\gamma=R$, $P=Q_\theta=R^{-1}h_{-it}(i,\pi)$ to the last logarithm in (\ref{eq:rotate_kernel}), then it is 
$$
-\log\left(\dfrac{i(t+1)\cdot \cos(\theta/2)+\sin(\theta/2)}{i(t-1)\cdot \cos(\theta/2)+\sin(\theta/2)}\right).
$$
We have $\varphi(t,\eta,\theta) = \Re\Phi_{\eta}(R^{-1}(iy(t,\eta)), Q_\theta)$; also, $\varphi(t,\eta,\pi) = \Re\Phi_{\eta}(iy(t,\eta), h_{-it}(i,\pi))$. To conclude the calculation for the second assertion of our Lemma, it remains to put $z = R^{-1}(iy(t,\eta))$ to (\ref{eq:rotate_kernel}).
$\blacksquare$

\medskip

\noindent {\bf Remark.} Since $h_{-it}(i,\pi)=(it,1)$, the second assertion of Lemma \ref{lemma:diffeo} implies, by a calculation, that 
\begin{equation}
	\label{eq:phitt}
	\varphi(t,t,\theta) = \dfrac{1}{2}\log\left({\dfrac{(1+\cos\theta)\cdot(t^2+2t)+2}{(1+\cos\theta)\cdot(t^2-2t)+2}}\right)+\log\left(\dfrac{2-t}{2+t}\right) -\dfrac{4-2t^2}{4t-t^3}.
\end{equation}

\section{Kernel $L_t(\cdot,\cdot)$ gives a semiclassical PDO}

\label{sec:pdo}

In this Section we reduce weighted quadratic means of $u^\CC$ on $\Sigma_t$ to a quadratic form given by a pseudodifferential operator and evaluated on $u$.

Take some  $t_1, t_2<1$ positive and close enough one to another, $t_1<t_2$; 
take $t\in (t_1, t_2)$. Pick $g\colon \R\to \R^+$ smooth, nonnegative and supported by $[t_1, t_2]$. Take any $a\in C_0^\infty(\Sigma_t)$ with $\supp a$ small enough. Consider operator with kernel 
\begin{multline*}
	L_t(z_1, z_2) := \int\limits_{t_1}^{t_2} d\eta\, g(\eta) \int\limits_{\Sigma_t} dS_t(P)\, \left(\dfrac{z_1-\tilde Z(P)}{\bar z_1 - Z(P)}\right)^\tau e^{-\tau c_\eta \cdot \cosh\dist(z_1, P)} a(P) \times \\ 
	\times e^{-\tau c_\eta \cdot \overline{\cosh\dist(z_2, P)}}\left(\dfrac{\bar z_2-\bar{\tilde Z}(P)}{z_2 - \bar Z(P)}\right)^\tau e^{-2\tau \varphi(t, \eta, \theta(P))} =\\=
	\int\limits_{t_1}^{t_2} d\eta\, g(\eta) \int\limits_{\Sigma_t}\, dS_t(P) K_\eta^\tau(z_1, P) a(P) \overline{K_\eta^\tau(z_2, P)} \cdot  e^{-2\tau\varphi(t, \eta,\theta(P))}, ~~~ z_1, z_2\in\HH.
\end{multline*}
Here, for $P\in\mathcal G_{1}\setminus \HH$, we write $\theta(P)$ for angular coordinate of $P$ in horocycle coordinates $(x,y,t,\theta)$; recall also that $dS_t(\cdot)$ is invariant Liouville measure $d\mu_L$ transferred from $\{H_1=1/2\}$ to $\Sigma_t$ by horocycle parametrization from Section \ref{sec:coords_flows}, function $\varphi$ has been defined in Lemma~\ref{lemma:max_point_exists}. We assume that $\supp g$ is small enough such that $\mathcal T_{z,t}$ from Lemma \ref{lemma:diffeo} is a diffeomorphism of $\supp g\times (\R\mmod2\pi)$ onto some closed neighborhood of $T_z\HH \cap \{H_{-1}=1/2\}$. 

\medskip


\medskip

Recall that $\mathcal S(\cdot,\cdot,\cdot)$ is defined in Lemma \ref{lemma:SelbergTransform}. 
From the second assertion of that Lemma we derive the following

\begin{predl} If $D^\tau u=s^2u$ then 
	\label{predl:kernel_Selberg_mult}
	\begin{multline*}
	\int_{\HH}\int_{\HH}u(z_1) \bar u(z_2) L_t(z_1, z_2)\,d\mathcal A_2(z_1)\,d\mathcal A_2(z_2) =\\= \int_{\Sigma_t} |u^\CC(P)|^2 a(P)
	\left(\int_{t_1}^{t_2}d\eta \, g(\eta) |\mathcal S(\eta,\tau,s)|^2\cdot e^{-2\tau\varphi(t, \eta, \theta(P))}\right) \,dS_t(P).
	\end{multline*}
\end{predl}

We thus put 
\begin{equation}
\label{eq:B_define}
	B(P) := \int_{t_1}^{t_2}d\eta \,g(\eta) |\mathcal S(\eta,\tau,s)|^2\cdot e^{-2\tau\varphi(t, \eta, \theta(P))}
\end{equation}
for $P\in \Sigma_t$ so that $B^{-1/2}$ will govern the asymptotics of $u^\CC$.

\medskip

\noindent {\bf Remark.} We may right now notice that $B(P)\neq 0$  if $g\not\equiv 0$ is non-negative. Indeed, otherwise $\mathcal S(\eta,\tau,s)$ vanishes at a non-degenerate interval of $\eta$'s. Using \cite{Fay} we have 
\begin{equation}
	\label{eq:S_explicit}
	\begin{aligned}
&\mathcal S(\eta,\tau,s) = \int_1^{+\infty}  e^{-\tau c_{\eta}\cdot \cosh r} \mathcal P_{s,\tau}(r)\,d\cosh r, \\
&\mathcal P_{s,\tau}(r)=(1-\tanh^2r/2)^{\tilde s}\cdot{}_2F_1(\tilde s-\tau, \tilde s+\tau,1;\tanh^2r/2), ~~ \tilde s(\tilde s-1)=-s^2.
\end{aligned}
\end{equation}
For $\tilde s, \tau$ fixed, we have $\mathcal P_{s,\tau}(r) = O(e^{Nr})$ as $r$ approaches $+\infty$ for some $N$ large enough (see \cite[\S15.4(ii)]{DLMF}).  We then conclude that $\mathcal S(\eta,\tau,s)$ is analytic in $c_\eta$ and thus has at most a discrete set of zeroes. An asymptotics for $B$ will be derived in Proposition~\ref{predl:B_asymp} below.
\medskip


The following Proposition is our main assertion relating distribution of $|u^\CC|^2$ at $\mathcal G_{1}$ to microlocal distribution of $u$ at $T^*\HH$. We state it in a general form forgetting that $u$ is an eigenfunction for $D^\tau$ and that the sequence $\{u\}$ is quantum ergodic.

\begin{predl}
	\label{prop:pdo}
Let $u=u_n\colon \HH\to\CC$, $n=1,2\dots$, be functions such that
\begin{enumerate}
	\item $\sup\limits_{n\in\mathbb N, \, z\in \HH} \|u_n\|_{L^1(\mathcal B_\HH(z,1))}<+\infty$,
	\item for any compact $\mathcal K\subset \HH$, $\sup\limits_{n\in\mathbb N} \|u_n\|_{L^2(\mathcal{K})}<+\infty$.
\end{enumerate} 

There exists a smooth 
function $b_{1,t}(z,\xi_1,\xi_2)\in C^\infty(T^*\HH)$ depending smoothly also on $t$ but not depending on $u_n$ and $g$ with the following property:

Let $a\in C_0^\infty(\Sigma_t)$ be smooth with support small enough 
whereas $g\in C_0^\infty(\R)$ be supported by $[t_1, t_2]$ with $t_1<t<t_2$ and $t_1, t_2$ close enough to $t$. Put 
$$
\mathscr s(z,\xi_1, \xi_2) := b_{1,t}(z, \xi_1, \xi_2)\cdot g(\Eta_{z,t}(-\xi_1, -\xi_2))\cdot 
a(Q(z,t,\Eta_{z,t}(-\xi_1, -\xi_2),\Theta_{z,t}(-\xi_1, -\xi_2)))
$$ 
with $Q$ being defined at Lemma \ref{lemma:max_point_exists}. Let $\tau=\tau_n\to+\infty$,  $A := \Op_{1/\tau_n}\mathscr s$ be semiclassical PDO with symbol $\mathscr s$. Then
\begin{equation}
	\label{eq:quadratic_form_equiv}
		\int_{\HH}\int_{\HH}u_n(z_1) \bar u_n(z_2) L_t(z_1, z_2)\,d\mathcal A_2(z_1)\,d\mathcal A_2(z_2)=O(1/\tau_n^4)+
	\tau_n^{-3}\cdot \langle Au_n,u_n \rangle_{L^2(\HH)}.
\end{equation}
\end{predl}

Notice that $\Eta_{z,t}(-\xi_1, -\xi_2)$, $\Theta_{z,t}(-\xi_1, -\xi_2)$ are initially defined for $(z, \xi_1, \xi_2)$ near $\{H_1=1/2\}$. By making $\supp g$ sufficiently small we can make symbol $\mathscr s$ well defined for all $\xi_1$, $\xi_2$ via continuation by zero.

The proof of Proposition \ref{prop:pdo} is obtained by standard tools giving Composition Theorem in the theory of Fourier Integral Operators. This theory does not cover our case of semiclassical operators with complex phase. If the corresponding canonical graphs calculus is established then we will be able to argue as follows. Let $A_1\colon (\mbox{functions on }\HH)\to(\mbox{functions on }\Sigma_t)$ be operator with kernel $K_t^\tau(z,P)\cdot \exp(-\tau\varphi(t,t,\theta(P)))$ , and $\mathcal M_a\colon (\mbox{functions on }\Sigma_t)\to(\mbox{functions on }\Sigma_t)$ be multiplication by $a$. Then, in the left-hand side of relation claimed in Proposition~\ref{prop:pdo} we have quadratic form given by $A_1^*\mathcal M_a A_1$, up to $\eta$-mollification as above; this mollifier is indeed necessary, for, otherwise, there will be dimension defect in graphs. More precisely, operator with kernel 
$$
\tilde L(z_1, z_2) := \int\limits_{\Sigma_t} dS_t(P)\, K_t^\tau(z_1, P) a(P) \overline{K_t^\tau(z_2, P)} \cdot  e^{-2\tau\varphi(t, t,\theta(P))}
$$ will not be a semiclassical PDO since its symbol is too singular.

Instead of our way of smoothing, we may put mollifier in $A_1$ replacing the latter operator with  operator 
$A_2$ having kernel 
$$
\int_{t_1}^{t_2}d\eta\, g(\eta) K_\eta^\tau(z,P)\cdot e^{-\tau\varphi(t,\eta,\theta(P))}, ~~ z\in\HH, ~ P\in\mathcal G_1.
$$
Replace $A_1^*\mathcal M_aA_1$ by $A_2^*\mathcal M_aA_2$, this operator also would allow us to study eigenfunctions. Now  we wish to be able to consider all the three factors as Fourier Integral Operators --- semiclassical and with complex phases. FIO's Composition Theorem (still unproved for such operators)  would allow to compose their graphs (see Introduction) and arrive to the identical graph for the whole $A_2^*\mathcal M_a A_2$. Together with symbol multiplication, this would lead us to a semiclassical PDO at the right-hand side of the relation from Proposition~\ref{prop:pdo}. 

Let us compare our approach to arguments from \cite{Ze07} in more details. Let $M$ be real-analytic Riemann manifold. In \cite{Ze07},  Zelditch deals with 
eigenfunctions $u^{\geod}_\lambda$ with $-\Delta u_\lambda^{\geod}=\lambda u^{\geod}_\lambda$, 
$\lambda\to+\infty$ ranges the set of 
eigenvalues of $-\Delta_M$, minus Laplacian on~$M$; here, $\hbar=1/\sqrt{\lambda}\to 0$ is the typical wavelength. 
Such functions admit analytic continuation into  (geodesic) Grauert tube $\mathcal G^{\geod}_t$ with radius $t>0$ small enough, denote the latter continuations by $u^{\geod, \CC}_\lambda$. To study them, we have to work with half-heat semigroup  $\exp\left(-t\sqrt{-\Delta_M}\right)$; it can be considered as imaginary time evaluation of half-Schr\"odinger evolution $\exp\left(it\sqrt{-\Delta_M}\right)$. The latter is known to have canonical graph given by the graph of geodesic flow. We pass to complex time. Let $t$ be small. In $L^2(\partial\mathcal G_t^{\geod})$, consider the subspace $\mathcal O(\partial\mathcal G_t^{\geod})$ of CR-holomorphic functions. Then $\exp(-t\sqrt{-\Delta_M})$ can be understood as operator $L^2(M)\to \mathcal O(\partial\mathcal G_t^{\geod})$, by analytic continuation of its Schwartz kernel. Moreover, this is Fourier Integral Operator with complex phase.  
The graph of the latter operator is given by complexified geodesic flow. Consider, as above, for some $a\in C^\infty_0(\partial\mathcal G_t^{\geod})$, operator $\mathcal M_a$ acting on functions on $\partial\mathcal G_t^{\geod}$ via multiplication by $a$. Then we may apply FIO Composition Theorem for 
\begin{equation}
\label{eq:FIO_composition_Zelditch}
\exp(-t\sqrt{-\Delta_M})^*\mathcal M_a\exp(-t\sqrt{-\Delta_M})\colon L^2(M)\to L^2(M).
\end{equation}
Let $\mathscr B_r^*M$ denote co-ball bundle over $M$ consisting of covectors of Riemannian lengths less than some $r>0$ small enough. Let $h^{\geod}_{it}\colon 
\mathscr B^*_rM\to M^\CC$ be geodesic flow on $M$ with complex time and values in complexification $M^\CC$ of $M$. Graph composition for~(\ref{eq:FIO_composition_Zelditch}) leads to the identical graph; further, the operator as above, in the leading order, reduces to PDO with symbol $\tilde{\mathscr s}$ coinciding to  $a\circ h^{\geod}_{it}$ at unit spherical bundle over $M$. 
Namely, \emph{localized} behavior, as $\lambda\to\infty$,  of $u_\lambda^{\geod,\CC}$ which is asymptotics of
$$
\int\limits_{\partial\mathcal G_t^{\geod}} a\cdot|u^{\geod,\CC}_\lambda|^2\,d\mu_L
$$
(with $\mu_L$ being, say, natural Liouville measure) is reduced to asymptotics of 
$$
\left\langle(\Op_{\hbar}\tilde{\mathscr s})u_\lambda^{\geod},u_\lambda^{\geod}\right\rangle_{L^2(M)}
$$
with $\hbar=1/\sqrt{\lambda}$; 
the latter is \emph{microlocalization} of $u_\lambda^{\geod}$. Since $u_\lambda^{\geod,\CC}$ is proportional to $\exp\left(-t\sqrt{-\Delta_M}\right)u_\lambda^{\geod}$, we thus may study growth 
of $u_\lambda^{\geod,\CC}$ in geodesic Grauert tube: the leading term in asymptotics for $|u^{\geod,\CC}_\lambda|$ on $\partial\mathcal G_t^{\geod}$, as $\lambda\to+\infty$, is $\exp(t\sqrt\lambda)$. Here, $t$ is instead of our answer $-B_0/2$. The role of function $t$ in geodesic setting has already been discussed at Section \ref{sec:coords_flows}. Notice also that Zelditch's asymptotic estimates as well as ours are true slice-wise, that is, at each $\partial\mathcal G^{\geod}_t$.

Other natural operators in this context are orthogonal Szeg\H{o} projector $L^2(\partial\mathcal G_t^{\geod})\to \mathcal O(\partial\mathcal G_t^{\geod})$ and derivative on functions on $\partial\mathcal G_t^{\geod}$ along Hamiltonian (Reeb) vector field. See  recent papers \cite{CR21},\cite{CR22} 
on Schwartz kernel asymptotics  of Toeplitz truncation of the latter vector fields. 

Operator in Zelditch's approach resembles operator $A_2^*\mathcal M_aA_2$ mentioned above rather than our $A_1^*\mathcal M_aA_1$ mollified. Our operators are not obtained by a matrix exponential. They are just something feasible to calculate. First, we need certain first-order relations making $h_{-it}(z,\theta)$ at least a stationary point of $|K_t(z,\cdot)|$ on $\Sigma_{t,\theta}$, and such first-order relations were provided by the appropriate choice of $c_t$. Second, we need $h_{-it}(z,\theta)$ to be the global maximum point of kernel absolute value. Third, we need summability conditions providing the first assertion of Lemma \ref{lemma:SelbergTransform}. All the operators having kernel of the form (\ref{eq:general_kernel}) possess rotational automorphy; we just deal with the simplest of them the and this leads to success.

What we lack, comparing to \cite{Ze07}, is homogeneity of symbol of operator $D^\tau$. Working with Laplacian, we just scale the same operator as $\hbar\to 0$. In our paper, we thus need a semiclassical family of operators, and they have complex phases. Since there is no Composition Theorem for our case, we apply perturbed complex stationary phase directly. 
To start, we outline the scheme of the argument. A detailed proof with more technicalities is given at the Appendix (Subsection \ref{subsec:Section4proofs}).

\medskip

\noindent {\bf Scheme of the proof of Proposition \ref{prop:pdo}.} Kernel $K_\eta^\tau(z,P)\cdot e^{-\tau\varphi(t,\eta,\theta(P))}$ does not exceed $1$ in absolute value; thus we expect that the left-hand side of (\ref{eq:quadratic_form_equiv}) is  constant-scale, up to degrees of $\tau$. 

By a direct  estimation one can see that  the contribution to the quadratic form of $(z_1,z_2)$  with $z_1$ or $z_2$ far enough from $\supp a$ is small for $\tau$ large. We proved similar estimates in the first assertion of Lemma \ref{lemma:SelbergTransform} assuming that $\sup|u|$ is finite. This is not assumed to be uniform even along a sequence of quantum ergodic functions $\{u_n\}$. But an examination of the proof of Lemma \ref{lemma:SelbergTransform} shows that  if $\sup\limits_{n\in\mathbb N, \, z\in \HH} \|u_n\|_{L^1(\mathcal B_\HH(z,1))}<+\infty$ then $z_1$ or $z_2$ far enough from $\supp a$ affect the left-hand side of (\ref{eq:quadratic_form_equiv}) by a negligibly small value, less than $\exp(-C\tau)$ with any given $C<\infty$.  

The same concerns the case when $z_1$ is separated from  $z_2$ (Lemma \ref{lemma:kernel_define_and_max}, in fact, provides strict non-degenerate maximum of $|K_t|$ at the given point).

We thus may assume that $z_1$ is close enough to $z_2$ and both range a compact set. Next step is to slice $\Sigma_t$ in integral for $L_t$ into $\bigcup\limits_{\theta\in\R\mmod 2\pi}\Sigma_{t, \theta}$, see Definition \ref{def:Grauert_sets}. Two-dimensional set $\Sigma_{t,\theta}$ is endowed with measure $\mathcal A_{2,t}$ which is the push-forward of hyperbolic area $\mathcal A_2$ under parametrization $\HH\ni z\mapsto h_{-it}(z,\theta)\in \Sigma_{t,\theta}$. We have estimate $|K_\eta(z_j, P)e^{-\varphi(t, \eta, \theta)}| \le 1$ for $P\in\Sigma_{t, \theta}$, $j=1,2$, turning to the equality at $P=Q(z_j,t,\eta,\theta)$. This point is $h_{-it}(z_j,\theta)$ if $\eta=t$. 
We are able  to apply perturbational complex stationary phase method as stated in \cite{Treves2} to 
\begin{equation}
\label{eq:Lteta}
\int\limits_{\Sigma_{t,\theta}} d\mathcal A_{2,t}(P) \, K_\eta^\tau(z_1, P)a(P)\overline{K_\eta^\tau(z_2, P)} \cdot e^{-2\tau\varphi(t, \eta, \theta)}
\end{equation}
with $\eta, \theta$ fixed. That is, we make use of almost-analytic continuations of amplitude and phase. Considering $z_1$ as parameter we complexify integration domain in (\ref{eq:Lteta}). If we start with $z_1=z_2$ then $P=Q(z_2,t,\eta,\theta)$ is stationary point in (\ref{eq:Lteta}). When we perturb this $z_1$,  the latter stationary point 
moves to the complexification of integration domain. This leads to the asymptotics of the form
$$
1/\tau\cdot\mathscr a_1(z_1, z_2, t, \eta, \theta)e^{\tau\cdot\Psi(z_1, z_2, t,\eta,\theta)}
$$
for (\ref{eq:Lteta}) with some $\mathscr a_1$ and $\Psi$; to be perfect, we need more terms of asymptotics --- up to $O(1/\tau^4)$ remainder. For phase obtained, $\left.\Psi\right|_{z_1=z_2}=0$. Also, a calculation shows that $$\left.d_{z_1}\right|_{z_1=z_2}{\Psi(z_1, z_2, t,\eta,\theta)}=i\mathcal T_{z_2,t}(\eta, \theta).$$
For estimation, it is useful to notice that $\Re\Psi\le 0$ due to \cite[Lemma X.2.5]{Treves2}.

Now replace $L_t(z_1, z_2)$ by $I:=1/\tau\cdot \int_{\R} d\eta\, g(\eta)\int_0^{2\pi }d\theta \, \mathscr a_1(z_1, z_2, t, \eta, \theta)e^{\tau\cdot\Psi(z_1, z_2, t,\eta,\theta)}$.  Localization in $I$ and repeated integration by parts show that contribution of $(z_1,z_2)$ with $|z_1-z_2|\ge \tau^{-2/3}$ to
$$
\int_{\HH}\int_{\HH}I\cdot u(z_1)\bar u(z_2)\,d\mathcal A_2(z_1)\,d\mathcal A_2(z_2)
$$ 
is $O(\tau^{-N})$ for any given $N<\infty$; to find an appropriate direction of this integration by parts, we may apply non-degenerateness provided by Lemma \ref{lemma:diffeo}.  

Now, assume that  $|z_1-z_2|\le \tau^{-2/3}$. In $I$, using Lemma \ref{lemma:diffeo} again, change variables as $(\eta,\theta)\mapsto \mathcal T_{z_2,t}(\eta, \theta)$; for this, we have, of course, to assume that $\eta$ is close enough to~$t$. Take long enough Taylor expansions over $z_1-z_2$ for phase and amplitude.  Principal term 
$$
1/\tau\cdot \int_{\R} d\eta\, g(\eta)\int_0^{2\pi }d\theta \,\mathscr a_1(z_2, z_2, t, \eta, \theta)e^{\tau\cdot i\mathcal T_{z_2}(\eta,\theta)[z_1-z_2]}
$$
leads to PDO from the statement. (Square brackets mean application of a covector to a vector.) All the other terms are negligible in the sense of quadratic forms by 
Calderon--Vailliancourt Theorem. Finally, sign before $\xi_1$ and $\xi_2$ appears during examination of quadratic form given by the reduced kernel. 
$\blacksquare$


\medskip

Now assume  that functions $u$, as in Theorem\ref{th:growth}, are uniformly distributed at 
$\{H_1=1/2\}$. This implies that, for $A$ as in Proposition \ref{prop:pdo} and for $Q=Q(z,t,\Eta_{z,t}(-\xi_1, -\xi_2),\Theta_{z,t}(-\xi_1, -\xi_2))$, we have 
\begin{multline}
	\label{eq:limit_pass_que}
\hspace{-0.1cm}
\langle Au,u \rangle_{L^2(\HH)} \to \int\limits_{\{H_1=1/2\}} b_{1,t}(z, \xi_1, \xi_2)\cdot g(\Eta_{z,t}(-\xi_1, -\xi_2))\cdot a(Q) 
\,d\mu_L(z,\xi_1, \xi_2)
=\\=
g(t)\cdot\int_{\{H_1=1/2\}} b_{1,t}(z, \xi_1, \xi_2)\cdot a(h_{-it}(z,\Theta_{z,t}(-\xi_1,-\xi_2)))\, d\mu_L(z,\xi_1, \xi_2)
\end{multline}
as $\tau\to+\infty$ (see (\ref{eq:EtaThetaCritical})). 
Recall that $\mu_L$ is appropriately normed Liouville measure on $\{H_1=1/2\}$, see Section \ref{sec:coords_flows}. Since we may take arbitrary $a$, we conclude that $\tau^3\cdot |u^\CC(P)|^2\cdot B(P) \cdot dS_t(P)$ converge to a measure mutually absolutely continuous with respect to $dS_t(P)$. In other words, as functions $u$ become equidistributed at $\{H_1=1/2\}$, functions $\tau^3\cdot |u^\CC(P)|^2\cdot B(P)$ become equidistributed on $\Sigma_t$ --- up to a smooth non-vanishing factor.

Convergence in (\ref{eq:limit_pass_que}) is uniform when $a$ ranges some compact set of symbols, namely, when all the derivatives of $a$ up to some sufficient order are bounded. Also, this convergence is uniform when $t$ ranges a compact subset in $(0,1)$. The same concerns limit relation from  Proposition \ref{prop:pdo}. Thus, integration of result of that Proposition over~$t$  leads us to the following

\begin{predl}
	\label{prop:weakstarconvergence}
	Let $0<t_1<t_2<1$ with $t_1$ close enough to $t_2$. 
	There exists a smooth strictly positive function $b_2\colon \mathcal G_{t_2}\setminus \clos\mathcal G_{t_1}\to \R^+$ with the following property: 
	
	Let $g\in C_0^\infty([t_1, t_2])$ and $B$ be as defined in \emph{(\ref{eq:B_define})}. 
	For such $B$, for a sequence of quantum ergodic functions $\{u\}$ as in Theorem \ref{th:growth} and for their complexifications $\{u^\CC\}$ we have
	$$
	\tau^3 |u^\CC(P)|^2 \cdot B(P) \xrightharpoondown[\tau\to+\infty]{}^* g(t(P))\cdot b_2(P).
	$$
	Here, weak* convergence is understood as in Theorem \ref{th:growth}.
\end{predl}

\noindent {\bf Remark.} In the following Section, we will ensure that $B(P)/g(t(P))$ asymptotically  does not depend on $g$ as $\tau\to+\infty$, this is natural to expect. 

\section{Asymptotics for $B$}

\label{sec:B_asymp}

Now we calculate asymptotics for $B$ when $\tau$ is large. Notice, by the way, that this is not necessary to prove Theorem \ref{th:main}. Since we are going to apply Lelong--Poincar\'e formula to arrive to that Theorem, we may just prove that $\dfrac{2\log|u^\CC|}{\tau} +\dfrac{\log B}{\tau}\xrightarrow{\tau\to\infty}0$ (see Lemma~\ref{lemma:psh} below; the argument can be modified for rather implicit $B$). Then 
it remains to find asymptotics for $\dfrac{\partial \log B(P)}{\partial t(P)}$, $\dfrac{\partial\log B(P)}{\partial \theta(P)}$.  
If $g\ge 0$ then this can be done by differentiating (\ref{eq:B_define}) or rather only the exponential function therein since the integrand is non-negative in this case. The asymptotics of the quotient does not depend on the choice of $g$ and is clear if $\supp g$ tends to one-point set $\{t\}$. Knowledge of $d\log B$ is enough to apply Lelong--Poincar\'e formula.

To calculate asymptotics for $B$, we start with asymptotic expression for $\mathcal S(\eta,\tau,s)$. To this end, formulae (\ref{eq:S_explicit}) seem to be unuseful. Indeed, we may try to represent hypergeometric function therein by an integral expression in the spirit of \cite[\S15.6]{DLMF}, then we have double integral for~$\mathcal S$. The maximum of exponential expression therein seems to be always on the boundary of $2$-dimensional contour of integration; also, this maximum \emph{does not} lead to the correct answer which is strictly less: boundary asymptotics should necessarily cancel, and this cannot be eliminated by a deformation of the contour. 

Instead, we make use of the spectral nature of $\mathcal S$ and of geometric intuition elaborated by now:

\begin{predl}
	\label{predl:S_asymp}
	As $\tau\to+\infty$ and $s=o(\tau)$, we have 
	\begin{equation*}
		|\mathcal S(\eta,\tau,s)| \sim \tau^{-1} \cdot b_3(\eta) \cdot \exp \tau\varphi(\eta,\eta,\pi) 
	\end{equation*}
with some $b_3$ smooth and separated from zero for $\eta$ strictly inside of $(0,1)$. The quotient of left- and right-hand sides tends uniformly to $1$ for such $\eta$. 
\end{predl}

\noindent {\bf Proof.} 
Notice that if $s_1= \sqrt{s^2-1/4}$, $v(z) = \left(\Im z\right)^{\frac12+is_1}$ ($z\in\CC$) then $D^\tau v = s^2 v$. It is possible to check that Lemma \ref{lemma:SelbergTransform} is still valid for such $v$. Then, since $v(i)=1$, we have
$$
\mathcal S(\eta,\tau,s)  = \int_{\HH} v(z) K^\tau_\eta(z,i)\,d\mathcal A_2(z).
$$
Put $s_2 := \bar s_1$ (we do not assume $s>1/4$). Observe that $K^\tau_\eta(z,i)=\overline{K^\tau_\eta(i,z)}$ for $z\in\HH$. Therefore
$$
\bar{\mathcal S}(\eta,\tau,s) = \int_{\HH}y^{1/2-is_2}K_\eta^\tau(i,z)\,d\mathcal A_2(z).
$$

The integrand admits a continuation to $\mathcal G_1$ analytic with respect to complexified components of $z$. Thus, we are able to make use of high-dimensional steepest descent method as stated in \cite{Fedoryuk}: we are going to shift a contour of integration having real dimension $2$ in $2$-dimensional complex space in order to hit a stationary point.

The integral for $\bar{\mathcal S}(\eta,\tau,s)$ is
$$
\int_{\HH} Y(P)^{-\frac32-is_2} \left(\dfrac{i-\tilde Z(P)}{-i - Z(P)}\right)^\tau\cdot \exp(-\tau c_\eta \cosh\dist(i,P)) \, dX(P)\wedge dY(P).
$$
Here $P$ ranges $\HH$ but we may consider $P\in \mathcal G_1$. Denote by $\omega$ the analytic $2$-form under integral sign. By (\ref{eq:ReY_ImX}), $\omega$ is well-defined at $\mathcal G_1$. 

We have $h_{-it}(x+iy,\pi) = (x+ity,y)$, and $t$ here can be understood as a homotopy parameter. Pick $r>0$ large enough. Consider Euclidean rectangle $\mathcal R:=[-r^{100},r^{100}]\times [1/r,r]\subset \HH$. 
In $\mathcal G_1$, consider contour $\mathscr M$ of real dimension $2$ consisting of the following parts:
\begin{gather*}
\HH\setminus \mathcal{R};\\
\mathcal R_\eta := h_{-i\eta}(\mathcal R,\pi)
\subset \Sigma_{\eta,\pi};\\
M_1 := \left\{h_{-it}(x+iy,\pi)\colon y\in\{r,1/r\}, \, x\in[-r^{100},r^{100}], \, t\in[0,\eta]\right\}\subset \mathcal G_1;\\
M_2 := \left\{h_{-it}(x+iy,\pi)\colon x=\pm r^{100}, \, y\in[1/r,r],\, t\in[0,\eta]\right\}\subset \mathcal G_1.\\
\end{gather*}
Under appropriate orientation of these parts, $\mathscr M$ is homotopic to $\HH$ with $\HH\setminus\mathcal R$ fixed during the homotopy. Therefore, $\int_{\HH}\omega = \int_{\mathscr M}\omega$.

We claim that, for any given $C>0$, there exists $r$ large enough such that $\int_{(\HH\setminus\mathcal R)\cup M_1\cup M_2}\omega = O(e^{-C\tau})$ as $\tau\to+\infty$.

To begin, notice that $\int_{M_1}\omega=0$ since $dX\wedge dY=0$ on $M_1$. (Here, estimation differs from that of the proof of Lemma \ref{lemma:SelbergTransform}, case $|x|< C$, $y<\epsilon$.) 
Further, $\int_{M_2}\omega$ is a sum of two expressions of the form 
\begin{equation*}
	\pm i\int\limits_{1/r}^rdy\int\limits_0^\eta dt\, y^{-1/2-is_2}
	\left(\frac{-x+i(1+y-ty)}{-x-i(1+y+ty)}\right)^\tau
	\cdot\exp\left(-\tau c_\eta -\tau c_\eta\cdot\frac{(x+ity)^2+(1-y)^2}{2y}\right).
\end{equation*}
with $x=\pm r^{100}$. The first factor is $\le r^{}$ in absolute value, the second one has modulus $\le 1$ (recall that $\tau>0$); the expression under exponential in the third one is $\le -\tau\cdot c_\eta r^{197}$ for $r$ large enough. This leads to the desired.

To manage with $\int_{\HH\setminus\mathcal R}\omega$, apply polar coordinates in $\HH$. Gauge factor is unimodular in this case. Let $z\in\HH$ be integration variable, $\rho\ge 0$ be $\dist(i,z)$. Then $\Im z\le e^{\rho}$. Using polar coordinates $(\rho,\theta)$ with metric tensor $d\rho^2+\sinh^2\rho\,d\theta^2$, we estimate 
$$
\left|\int_{\HH\setminus\mathcal R}\omega\right|\le 2\pi\int_{\rho_0}^\infty d\rho\, e^{\rho}\cdot \exp(-\tau c_\eta\cosh\rho)\cdot \sinh\rho
$$
where $\rho_0\in(0,+\infty)$ can be taken arbitrarily large by appropriate choice of $r$. The latter inequality implies that $\left|\int_{\HH\setminus\mathcal R}\omega\right|\le e^{-C\tau}$ if $r$ is large enough.

We thus conclude  that $\bar{\mathcal S}(\eta,\tau,s)=\int_{\mathcal R_\eta}\omega+O(e^{-C\tau})$ for $r$ large. By a direct calculation we check that $P_0=h_{-i\eta}(i,\pi)=(i\eta,1)$ is a stationary point of phase $\Phi_\eta(i,\cdot)=\log K_\eta(i,\cdot)$ at the whole $4$-dimensional complexified Lobachevsky plane\footnote{Surprisingly, here  we may not replace $\pi$ by an arbitrary $\theta\in\R\mmod2\pi$. Also, for $\theta=\pi$, we may prove stationarity for imaginary horocycle time $-i\eta$ by checking the same for real time $\eta$ instead and then by analytic continuation to imaginary time.}. 
This point is also non-degenerate: the determinant of 
$$
\begin{pmatrix}
	\dfrac{\partial^2 \Phi_\eta(P,i)}{\partial X(P)^2} & \dfrac{\partial^2 \Phi_\eta(P,i)}{\partial X(P)\partial Y(P)}\vspace{2mm} \\
	\dfrac{\partial^2 \Phi_\eta(P,i)}{\partial X(P)\partial Y(P)} &
	\dfrac{\partial^2 \Phi_\eta(P,i)}{\partial Y(P)^2} 
\end{pmatrix}
$$
is, by a calculation, $4\cdot\dfrac{4-3\eta^2}{\eta^6-8\eta^4+16\eta^2}$ which is non-zero. Finally, we observe that, by second assertion of Lemma \ref{lemma:kernel_define_and_max},  $P_0=\argmax{P\in\Sigma_{\eta,\pi}}{|K_\eta(i,P)|}$.

The last difficulty is that we have double asymptotics: besides $\tau$, there is also (possibly large) $s_2$ in our integral. But since $s=o(\tau)$, we may replace $Y(P)^{-\frac32-is_2}$ by $Y(P_0)^{-\frac32-is_2}$. Indeed, near $P_0$, deform surface of integration $\mathcal R_\eta\subset \Sigma_{\eta,\pi}$ to the canonical steepest descent contour $W$ as in the proof of \cite[Chapter V, \S1.3, Theorem 1.1]{Fedoryuk}; then the asymptotics by $\tau$ (with $s, s_2$ being fixed) is calculated by Laplace method. For any $\tau,s$ take $\epsilon>0$ such that $1/\tau=o(\epsilon)$, $\epsilon=o(1/s)$ as $\tau\to+\infty$, this is possible by the assumptions from Theorem~\ref{th:growth}. The contribution to the integral of $P\in W$ with $\dist(P,P_0)>\epsilon$ is negligible. (The distance is understood in some, say, Euclidean coordinates in $\HH^\CC$.) For the remaining part, recall that $Y(P_0)=1$. We have $\left|Y(P)^{-\frac32-is_2}-1\right|=O(s\epsilon)=o(1)$ as $\dist(P,P_0)\le \epsilon$. Since there is no oscillation in Laplace-type integral along canonical contour $W$, it is indeed safe to put $Y(P_0)(=1)$ instead of $Y(P)$ in $\int_W Y(P)^{-\frac32-is_2} K_\eta^\tau(P,i)\,dX\wedge dY$. 
This concludes the proof, the higher-dimensional saddle point method indeed leads to the proposed answer.
$\blacksquare$

\medskip

Now we compute the asymptotics for $B$. Remark after Proposition \ref{prop:weakstarconvergence} suggests that $\eta=t$ should be a Laplace point in integral (\ref{eq:B_define}) for $B$. This leads to the proof of the following 

\begin{predl}
	\label{predl:B_asymp}
	Let $t\in (0,1)$, and let $t_1, t_2\in(0,1)$ be close enough to $t$ and such that $t_1<t<t_2$. Let $\theta\in\R\mmod2\pi$, and let also $g\colon \R\to\R$ be smooth \emph{non-negative} function with support $[t_1,t_2]$, 
	$g(t)>0$ at $(t_1, t_2)$. 
	
	As $\tau\to +\infty$ and $s=o(\tau)$, for $P\in \Sigma_{t,\theta}$  we have 
	\begin{equation*}
	B(P) \sim \tau^{-5/2}\cdot b_4(t) g(t)\cdot \exp(\tau B_0).
	\end{equation*}
	Here, $b_4$ is some smooth function on $\mathcal G_1\setminus \HH$ separated from $0$. The quotient of left- and right-hand sides of this relation tends to $1$ uniformly by $t$ strictly inside of $(t_1,t_2)$. 
\end{predl}

\noindent {\bf Proof.} 
Notice that integrand in relation (\ref{eq:B_define}) defining $B$ is non-negative since $g$ is such. Thus we may put asymptotics obtained in Proposition \ref{predl:S_asymp} to (\ref{eq:B_define}). Using also the second assertion of Lemma \ref{lemma:diffeo}, we get
\begin{multline*}
B(P) \sim \tau^{-2}\cdot \int_\R g(\eta) b_3^2(\eta) e^{2\tau\varphi(\eta,\eta,\pi)-2\tau\varphi(t,\eta,\theta)}\,d\eta
=\\=
\tau^{-2}\cdot 
\exp(\tau B_0)\cdot \int_\R g(\eta) b_3^2(\eta) e^{2\tau\varphi(\eta,\eta,\pi)-2\tau\varphi(t,\eta,\pi)}\,d\eta.
\end{multline*}
To prove the required asymptotics, we thus need to show the following: if $t$ is fixed and 
$f(\eta)  = \varphi(\eta,\eta,\pi)-\varphi(t,\eta,\pi)$ then $f'(t)=0$, $f''(t)<0$. Indeed, then, if $\supp g$ is small enough then Laplace method leads to the desired.

Recall that $Q=Q(z,t,\eta,\theta)$ has been defined at Lemma~\ref{lemma:max_point_exists}. 
To take 
$$
\dfrac{\partial}{\partial\eta}\varphi(t,\eta,\pi) = \dfrac{\partial}{\partial\eta} \log|K_\eta(i,Q(i,t,\eta,\pi))|,
$$ 
notice that 
$$
\left.\dfrac{\partial}{\partial\eta_1}\right|_{\eta_1=\eta} \log|K_{\eta}(i,Q(i,t,\eta_1,\pi))|=0
$$ 
since $Q(i,t,\eta,\pi) = \argmax{P\in \Sigma_{t,\pi}} |K_{\eta}(i,P)|$ and $Q(i,t,\eta_1,\pi)\in \Sigma_{t,\pi}$ for any $\eta_1$ close enough to $\eta$. We thus find 
\begin{equation}
\label{eq:weird_diff}
\dfrac{\partial}{\partial\eta}\varphi(t,\eta,\pi)=
\left.\dfrac{\partial}{\partial\eta_1}\right|_{\eta_1=\eta} \log|K_{\eta_1}(i,Q(i,t,\eta,\pi))| = -\dfrac{dc_\eta}{d\eta}\cdot \Re\cosh\dist(i,Q(i,t,\eta,\pi)).
\end{equation}
If $\eta=t$ then we may proceed calculations using (\ref{eq:phitt}) and arrive to $f'(t)=0$.

To find $f''(t)$ we still use (\ref{eq:weird_diff}). For $\eta$ close enough to $t$ we may write $Q(i,t,\eta,\pi) = h_{-it}(x(\eta)+iy(\eta),\pi)$ with some $x(\eta)+iy(\eta)\in \HH$ depending smoothly on~$\eta$, $x(t)+iy(t)=i$. Calculating Hesse matrix of $(x,y)\mapsto \log|K_\eta(i,h_{-it}(x+iy,\pi))|$ we find $\dfrac{dx(\eta)}{d\eta}=0$, $\left.\dfrac{dy(\eta)}{d\eta}\right|_{\eta=t}=-\dfrac{t(4-3t^2)}{2(4-3t^2+t^4)}$, and, by (\ref{eq:weird_diff}), 
$$ 
f''(t) =\dfrac{3t^2-4}{t(t^4-3t^2+4)} <0. 
$$
This concludes our computational proof. $\blacksquare$

\medskip

Theorem \ref{th:growth} now is an immediate consequence of 
Proposition \ref{prop:weakstarconvergence} since both sides of the limit relation therein are non-negative. Moreover, we may give slice-wise version of Theorem \ref{th:growth}:

\begin{predl}
\label{predl:sliced_growth}
For $u_n, \tau_n$ as in Theorem \ref{th:growth} \emph{(}$\tau_n\to+\infty$\emph{)} and for each $t\in(0,1)$,
$$
\tau_n^{1/2} \cdot |u_n^\CC|^2 \cdot 
\exp(\tau_n B_0) \xrightharpoondown[\tau_n\to+\infty]{}^*  b \mbox{ in } \mathcal D'(\Sigma_t)
$$
with  smooth  $b>0$ defined at $\mathcal G_1\setminus\HH$ and also with the meaning of weak* convergence as in Theorem \ref{th:growth} \emph{(}but on $\Sigma_t$\emph{)}.
\end{predl}

\section{Logarithm of weak* convergence}
\label{sec:psh}

To derive Theorem \ref{th:main} from Theorem \ref{th:growth}, we have to take the logarithm of the result of the latter one. This is done by a rather standard trick with plurisubharmonic dichotomy.

Recall that $u^\CC=u_n^\CC$, $\tau=\tau_n$, $s=s_n$ depend on $n=1,2,\dots$.

\begin{lemma}
	\label{lemma:psh}
	We have 
	\begin{equation*}
		\dfrac{2}{\tau_n} \cdot \log|u_n^{\CC}|+ B_0 \xrightarrow{n\to\infty}0\mbox{ in } L^1_{\loc}(\mathcal G_1).
	\end{equation*}
\end{lemma}

\noindent {\bf Proof.} We mostly follow Zelditch (\cite{Ze07}).

By the definition of $\mathcal S$ (Lemma \ref{lemma:SelbergTransform}),  
\begin{equation*}
	\mathcal S(1/2, \tau_n, s_n)u_n^{\CC}(P) = \int_{\HH} u_n(z)K^{\tau_n}_{1/2}(z,P)\,d\mathcal A_2(z), ~~ P\in\mathcal G_1.
\end{equation*}
We may estimate the integral using the condition  $\sup\limits_{n\in\mathbb N, \, z\in \HH} \|u_n^{}\|_{L^1(\mathcal B_\HH(z,1))}<+\infty$ 
required in the statement of Theorem \ref{th:growth}. Using estimates from the proof of first assertion in Lemma \ref{lemma:SelbergTransform}, 
one is able to 
see that, for any compact set $\mathcal K\subset \mathcal G_{1}$,
\begin{equation*}
	\sup\limits_{n\in\mathbb N}\sup\limits_{P\in \mathcal K} \dfrac{\log |\mathcal S(1/2, \tau_n, s_n)u_n^{\CC}(P)|}{\tau_n} < +\infty
\end{equation*}
(see discussion in the proof of Proposition \ref{prop:pdo}). 
Since we already have  asymptotics for $\mathcal S(1/2,\tau_n, s_n)$ given by Proposition \ref{predl:S_asymp}, we may conclude that 
\begin{equation}
	\label{eq:u_upper}
	\sup\limits_{n\in\mathbb N}\sup\limits_{P\in \mathcal K} \dfrac{\log |u_n^{\CC}(P)|}{\tau_n}<+\infty,
\end{equation}
$\mathcal K$ being any fixed compact in $\mathcal G_1$.

Consider plurisubharmonic functions $\dfrac{\log|u_n^{\CC}|}{\tau_n}$, $n=1,2,\dots$. From (\ref{eq:u_upper}) we see that these functions are bounded from the above on any compact set in $\mathcal G_1$ uniformly by $n$. By \cite[Theorem 4.1.9]{Hor1}, we have the following plurisubharmonic dichotomy: either $\dfrac{\log|u_n^{\CC}|}{\tau_n}\xrightarrow{n\to\infty} -\infty$ uniformly on each compact subset in $\mathcal G_1$; or, up to subsequence of indices $n$, functions $\dfrac{\log|u_n^{\CC}|}{\tau_n}$ converge in $L^1_{\loc}(\mathcal G_1)$ as $n\to\infty$.

The first case is impossible. Indeed, this would contradict Theorem \ref{th:growth} since $B_0$ is bounded from the below on compacts in $\mathcal G_1$. 

We thus may suppose, up to subsequence, that $\dfrac{\log|u_n^{\CC}|}{\tau_n} \xrightarrow{n\to\infty} f$ in $L^1_{\loc}(\mathcal G_1)$ for some function $f\in L^1_{\loc}(\mathcal G_1)$. Let $f^*$ be upper-semicontinuous regularization of $f$ (\cite[Theorems 4.1.11, 4.1.8]{Hor1}). Then $f^*$ is plurisubharmonic and equals $f$ almost everywhere in $\mathcal G_1$ with respect to Euclidean coordinates therein. 

First, we are going to prove that $2f^*+B_0=0$ 
in $\mathcal G_1\setminus \HH$.


Let us show  that $2f^*+B_0 \le 0$ almost everywhere in $\mathcal G_1\setminus \HH$. Indeed, otherwise, passing to a subsequence converging almost everywhere and applying D.~Egorov Theorem, we may assume that $\lim\limits_{n\to\infty} \left(B_0+\frac1{\tau_n}\log|u_n^{\CC}|^2\right)$ exists, is uniform and is greater or equal than some $\delta>0$ on a measurable set $E\subset\mathcal G_1\setminus \HH$ of a positive measure. Then $|u_n^{\CC}|^2 \cdot \exp (\tau_nB_0) \ge \exp(\tau_n\delta/2)$ on $E$ for $n$ large. 
We then arrive to contradiction to the weak*  convergence from Theorem \ref{th:growth}.

Now  prove that $2f^*+B_0 \ge 0$ in $\mathcal G_1\setminus \HH$. Suppose that $2f^*(P_0)+B_0(P_0)<-\delta$ for some $\delta>0$ and for some $P_0\in\mathcal G_1\setminus \HH$. Then, since $f^*$ is upper-semicontinuous and $B_0$ is continuous, we have  $2f^*(P)+B_0(P)<-\delta$ for $P$ in some  neighborhood $U$ of $P_0$ precompact in $\mathcal G_1\setminus \HH$; we may assume the same even  for $P\in\clos U$. Then, due continuity of $B_0$ again and by \cite[Theorem~4.1.9(b)]{Hor1}, 
\begin{equation*}
	\varlimsup\limits_{n\to\infty}\sup_{\clos U}\left( \dfrac{\log|u_n^{\CC}|^2}\tau+B_0\right) \le \sup_{\clos U} (2f^*+B_0) < -\delta,
\end{equation*}
and $|u_n^{\CC}|^2 \cdot \exp (\tau_n B_0) < \exp(-\tau_n\delta)$ on $U$ for $n$ large enough.
This again contradicts  Theorem \ref{th:growth}. 

So, by now, from upper-semicontinuity and plurisubharmonicity of $f^*$ and by continuity of $B_0$, we have $f^*=-\dfrac{B_0}{2}$ in $\mathcal G_1\setminus \HH$. But the above considerations  do not provide any information on the behavior of functions $u_n^{\CC}$ near $\HH$. (In all the preceding arguments we had to assume that $t$ is separated from zero to get uniform estimates of reminders.) In Theorem \ref{th:main}, we do not cut $\HH$ from $\mathcal G_{1}$. 

\emph{Function $-B_0/2$ is plurisubharmonic on the whole $\mathcal G_1$.} Indeed, it is such near any point in $\mathcal G_1\setminus\HH$ since it coincides to $f^*$ therein. Also, $B_0$ is continuous everywhere in $\mathcal G_1$. Finally, denote by $\mathcal B_{\CC}(0,r)$ the disc in $\CC$ centered in $0$ and having radius $r>0$; let also $\mathcal H^2$ be area measure on $\CC$. If $P_0\in\HH\subset\HH\times\HH$, $v\in\CC\times\CC$ is a vector tangent to $\HH\times\HH$ in $P_0$, $r>0$ is small enough then 
$$
-B_0(P_0)/2=0\le -\dfrac{1}{\pi r^2}\int_{\mathcal B_{\CC}(0,r)} \frac{B_0(P_0+\alpha v)}2\, d\mathcal H^2(\alpha)
$$
because the integrand at the right-hand side is non-positive (see (\ref{eq:B0_define})). Thus, $-B_0$ is also plurisubharmonic in any point on $\HH$; by localization (\cite[Theorem 1.6.3]{HorSV}) we conclude that $-B_0$ is plurisubharmonic in $\mathcal G_1$.

Thus, both $f^*$ and $-\dfrac{B_0}{2}$ are (pluri)subharmonic on $\mathcal G_1$ and they coincide on a set $\mathcal G_1\setminus\HH$ having full Euclidean measure therein. Then they generate the same distribution. But a (pluri)subharmonic function is \emph{uniquely} defined by its distribution (\cite[4.1.8]{Hor1}), therefore $f^*=-B_0/2$ everywhere at the whole $\mathcal G_1$.
Proof of Lemma is complete.
$\blacksquare$

\medskip

Now, to derive Theorem \ref{th:main} from Lemma \ref{lemma:psh}, it remains to apply Lelong--Poincar\'e formula to the obtained weak* convergence.

\appendix

\section{Appendix: some technical proofs}

\label{section:appendix}


\subsection{Proof of Lemma \ref{lemma:SelbergTransform}}

\label{subsec:Sec3proofs}


Let's prove the first claim. We fix $\tau$ and $c_t$. By applying an isometry we may assume that $X(P)$ is close to $it$ and $Y(P)$ is close to $1$ (see (\ref{eq:gauge_isometries})). 
Any derivative of integrand over components $X(P)$ or $Y(P)$ can be written as
\begin{equation*}
\left(\dfrac{z-\tilde Z(P)}{\bar z - Z(P)}\right)^{\tau-k_1}\cdot \exp\left(-\tau c_t \dfrac{(x-X(P))^2+(y-Y(P))^2}{2 y Y(P)}\right) \cdot\frac{\mathscr P(x,y,X,Y)}{(\bar z-Z)^{k_2}(yY)^{k_3}}\cdot u(z),
\end{equation*}
where $z=x+iy$, $Z=X+iY$, $\tilde Z=X-iY$, $\mathscr P$ is a polynomial, $k_1,k_2, k_3\in \mathbb N\cup\{0\}$. 

Using (\ref{eq:ReY_ImX}), we see that $\Im(z-\tilde Z), |\Im(\bar z-Z)| >y$.
Thus, by  Cauchy--Bunyakovsky--Schwartz inequality we see that it is enough to prove that 
\begin{equation}
	\label{eq:kernel_convergence}
	\int_{\HH} (1+|x|^k+y^k+y^{-k})\exp\Re\left(-\tau c_t \dfrac{(x-X(P))^2+(y-Y(P))^2}{2 y Y(P)}\right) \, d\mathcal A_2(x+iy) <+\infty
\end{equation}
for any $k=0,1,2\dots$. We take some $C, C'>100$ large enough and $\epsilon$ small enough and subdivide integration domain  as follows:
\begin{enumerate}
	\item $|x|<C$, $y<\epsilon$. If $|X(P)-it|$ is small enough then $\Re(x-X(P))^2\ge -\frac{t^2+1}{2}$, $\Im(x-X(P))^2$ is bounded above. If $Y(P)$ is close to $1$ then $(y-Y(P))^2$ is close to $1$. Then $\Re\left((x-X(P))^2+(y-Y(P))^2\right)$ is positive and separated from zero whereas $\left|\Im\left((x-X(P))^2+(y-Y(P))^2\right)\right|$ is bounded from the above. Thus $\arg\left((x-X(P))^2+(y-Y(P))^2\right)$ is separated from $\pm\pi/2$. Then $\Re\dfrac{(x-X(P))^2+(y-Y(P))^2}{2 y Y(P)}\ge \dfrac{\const}{y}$ provided that $Y(P)$ is close to $1$. Integral~(\ref{eq:kernel_convergence}) then converges over $\{|x|<C, y<\epsilon\}$.
	
	\item  $|x|>C$, $y<100$. We have $y-Y(P) = O(1)$, $\Re (x-X(P))^2 \ge \const\cdot |x|$, $|\arg (x-X(P))^2|$ can be forced to be arbitrarily small. Thus
	$$
	\Re\dfrac{(x-X(P))^2+(y-Y(P))^2}{2 y Y(P)}\ge  \const\cdot \dfrac{|x|}{y}.
	$$
	But 
	$$
	\int\limits_{x>C}dx\int\limits_{0}^{100} dy\, e^{-\frac{\const\cdot x}{y}}(x^k +y^{-k}) = \int\limits_{0}^{100} e^{-\frac{\const}{y}}\cdot\left(\mathcal P_1(y)+\mathcal P_2(1/y)\right)\, dy <+\infty,
	$$
	here $\mathcal P_1, \mathcal P_2$ are some polynomials.
	
	\item $|x|>C$, $100< y < C$. Then notice that $|\arg (x-X(P))^2|$, $|\arg (y-Y(P))^2|$ are separated from $\pi/2$. The rest is the same as in the previous case.
	
	\item $|x|>C$, $y>C$. Then both $\arg(x-X(P))^2$, $\arg(y-Y(P))^2$ can be forced to be arbitrarily small, and real parts of the expressions under $\arg$ are bounded from the below by $\const\cdot x^2$ and $\const \cdot y^2$, respectively. Then 
	$$
	\Re\dfrac{(x-X(P))^2+(y-Y(P))^2}{2 y Y(P)}\ge \const \cdot \dfrac{x^2+y^2}{y}.
	$$
	We have 
	\begin{multline*}
		\int\limits_{C}^\infty dx\int\limits_{C}^{+\infty} dy\, e^{-\const\cdot\frac{x^2+y^2}{y}}(x^k+y^k) \le \int_0^{\pi/2} d\phi \int_0^{+\infty} dr \, e^{-\frac{\const\cdot r}{\sin\phi}} r^{k+1}(\sin^k\phi+\cos^k\phi) \le \\ \le \pi\cdot\int_0^{+\infty} dr \, e^{-{\const\cdot r}}\cdot r^{k+1} <+\infty.
	\end{multline*}
	
	\item $|x|< C$, $y>C'$. Then $(x-X(P))^2 = O(C^2)$, $\Re (y-Y(P))^2\ge \const \cdot y^2$, $\arg(y-Y(P))^2$ can be forced to be arbitrarily small. We first pick $C$, then $C'$ large enough. The integral (\ref{eq:kernel_convergence}) is majorized by $C \cdot \int_{C'}^{+\infty} (\const+y^k) e^{-\const \cdot y}\, dy$ which is finite.
	
	\item $|x|<C$, $\epsilon<y< C'$. This a proper part of our integral and there is no convergence problem.
\end{enumerate}

The proof of the first assertion is thus complete. To prove the second one, one argues as in  \cite[Theorem 1.5]{Fay}. 
The third assertion of our Lemma follows from analyticity of $K_t(z,P)$ with respect to $X(P)$ and $Y(P)$. Proof of Lemma is complete.~$\blacksquare$

\subsection{Proof of Proposition \ref{prop:pdo}}

\label{subsec:Section4proofs}


We start our estimates assuming that $z_1, z_2$ in the left-hand side of (\ref{eq:quadratic_form_equiv}) range a compact set. For $\eta\in(t_1, t_2)$, $z_1, z_2\in \HH$, $\theta\in\R\mmod2\pi$, put
\begin{equation*}
	L_{t,\eta,\theta}(z_1, z_2):=\int\limits_{\Sigma_{t,\theta}} d\mathcal A_{2,t}(P) \, K_\eta^\tau(z_1, P)a(P)\overline{K_\eta^\tau(z_2, P)} \cdot e^{-2\tau\varphi(t, \eta, \theta)}
\end{equation*}	
where $d\mathcal A_{2,t}(h_{-it}(x+iy,\theta))=\dfrac{dx\,dy}{y^2}$ is hyperbolic area transferred to $\Sigma_{t,\theta}$.

First suppose that $z_2=z_1$. 
Then by Laplace method and by Lemma \ref{lemma:kernel_define_and_max}, we have:
$$
L_{t,\eta,\theta}(z_1, z_1) = b_1(t,\eta,\theta)\cdot \frac{1}{\tau}\cdot a(Q(z_1,t,\eta,\theta)) 
+O\left(\frac1{\tau^2}\right)
$$
with some $b_1$ smooth and separated from zero. Indeed, $$Q(z_1,t,\eta,\theta)=\argmax{P\in \Sigma_{t, \theta}} \, |K_\eta(z_1, P)|$$ and this stationary point is non-degenerate. 

Now, we are going to perturb this $z_2$ which is initially $z_1$; 
thus, we now assume that $z_2$ is close enough to $z_1$ whereas $\eta$ is close to $t$. By the definition of $\varphi$, we have  $|K_\eta(z_j, P)e^{-\varphi(t, \eta, \theta)}| \le 1$ for $P\in\Sigma_{t, \theta}$, $j=1,2$. We make use of almost analytic continuation technique from \cite{Treves2}. We thus write $\HH = \R\times (0,+\infty)$ and again consider it as a subset of $\CC\times\CC$. Parametrize $\Sigma_{t, \theta}$ as $\{h_{-it}(z,\theta)\colon z\in\HH\}$ and put $f(z) := h_{-it}(z,\theta)$. Let us suppose that $\supp a$ is close enough to $f(z_1)$. Functions $z\mapsto K_\eta(z_1, f(z))$, $z\mapsto \Phi_\eta(z_1, f(z))$ admit  holomorphic continuations from $z\in \R\times (0,+\infty)$ to some neighborhood of $f^{-1}(\supp a\cap \Sigma_{t, \theta})$ in $\CC\times\CC$. This is easily seen from the explicit horocycle parametrization of Grauert tube (in fact, linear by components of $z$) and from  (\ref{eq:half_kernel}); we keep the same notation $K_\eta$, $\Phi_\eta$ for these analytic continuations.

By \cite[Lemma X.2.3 and Remark X.2.1]{Treves2}, we may assume that $z\mapsto a(f(z))$ is \emph{almost-analytically} continued to $\CC\times\CC$ from $\HH$, denote this extension by $a_1=a_1(z)$, $z\in \CC\times\CC$. Also, denote by $\Phi_{(1)}(z_2, z, \eta)$ ($z$ ranges some neighborhood of $f^{-1}(\supp a\cap \Sigma_{t, \theta})$ in $\CC\times\CC$) an almost-analytic extension of $z\mapsto {\overline{\Phi_\eta(z_2,f(z))}}$ which also does exist by the same reason. 

Write, as before, $P\in \CC\times\CC$ as $P=(X,Y)$. By the third assertion of Lemma \ref{lemma:kernel_define_and_max}, for $z_2$ close to $z_1$ and $\eta$ close to $t$, there exists the unique $z_0=z_0(z_1, z_2, t, \eta, \theta)\in\CC\times\CC$ for which 
\begin{equation}
	\label{eq:stat_almost_analytic}
	\left.d_{(\Re X, \Re Y)}\right|_{P=z_0} \left(\Phi_\eta(z_1, f(P)) + \Phi_{(1)}(z_2, P,\eta)\right) = 0
\end{equation}
(two complex equations for two complex variables, not holomorphic but almost-holomorphic as $z\in\CC\times\CC$ approaches $\HH$).  This $z_0$ depends smoothly on its arguments, $z_0(z,z,t,\eta,\theta)\in\HH$ for any $z\in\HH$, and $h_{-it}(z_0(z,z,t,\eta,\theta),\theta)=Q(z,  t, \eta, \theta)$. 

Now, for the sake of further phase calculations, take 
\begin{equation}
	\label{eq:diff_phase1}
	\left.d_{z_1}\right|_{z_1=z_2} \left(\Phi_\eta(z_1, f(z_0(z_1, z_2, t, \eta,\theta))) + \Phi_{(1)}(z_2, z_0(z_1, z_2, t, \eta,\theta),\eta)\right).
\end{equation}
By (\ref{eq:stat_almost_analytic}) and by almost-analyticity  of all the functions, 
$$
\left.d_{z_0}\right|_{z_1=z_2} \left(\Phi_\eta(z_1, f(z_0(z_1, z_2, t, \eta,\theta))) + \Phi_{(1)}(z_2, z_0(z_1, z_2, t, \eta,\theta),\eta)\right)=0,
$$
we conclude that (\ref{eq:diff_phase1}) equals $\left.d_{z}\right|_{z=z_1} \Phi_\eta(z, Q(z_1, t, \eta, \theta))=i\mathcal T_{z_1,t}(\theta,\eta)$ (see Lemma~\ref{lemma:diffeo}).

By complex stationary phase method as stated in \cite[X.3]{Treves2} we have, as $\tau\to+\infty$,
\begin{multline}
	\label{eq:stat_phase}
	L_{t,\eta,\theta}(z_1, z_2)= O(1/\tau^4)+\\ 
	+\left(\dfrac1\tau\cdot a_1(z_0(z_1, z_2,t,\eta,\theta))\cdot \mathscr b_2(z_1, z_2, t, \eta,  \theta)
	+\dfrac{a_2(z_1, z_2,t,\eta,\theta)}{\tau^2} 
	+
	\dfrac{a_3(z_1, z_2,t,\eta,\theta)}{\tau^3}
	\right)\times\\ \times\exp\left(\tau\cdot\left(\Phi_\eta(z_1, f(z_0(z_1, z_2, t, \eta,\theta))) + \Phi_{(1)}(z_2, z_0(z_1, z_2, t, \eta,\theta),\eta)-2\varphi(t,\eta,\theta)\right)\right).
\end{multline}
Here, $\mathscr b_2$ 
is a smooth function of its arguments 
separated from zero; $a_2$ and $a_3$   are linear combinations with smooth coefficients of derivatives of $a_1=a_1(z)$ (the almost-analytic continuation of $a$) with respect to components of $z\in\CC\times\CC$ taken at $z=z_0(z_1,z_2,t,\eta,\theta)$. The constant in the remainder $O(1/\tau^4)$ in (\ref{eq:stat_phase}) depends only on suprema of derivatives of the functions $a$ and also $K$ and $\varphi$ up to some  finite order.

Now integrate the obtained expression by $d{\eta}$ and $d\theta$. We expect that if $z_1$ is close to $z_2$ then the main contribution to  $L_t(z_1, z_2) = \int_0^{2\pi}d\theta\, \int_{\R}^{} d\eta\, g(\eta) L_{t,\eta,\theta}(z_1, z_2)$ will be given by 
\begin{multline*}
	\dfrac1\tau\cdot a_1(z_0(z_1, z_2,t,\eta,\theta))\cdot \mathscr b_2(z_1, z_2, t, \eta,  \theta)\times\\
	\times\exp\left(\tau\cdot d_{z_1}|_{z_1=z_2}\left(\Phi_\eta(z_1, f(z_0(z_1, z_2, t, \eta,\theta))) + \Phi_{(1)}(z_2, z_0(z_1, z_2, t, \eta,\theta),\eta)\right)[z_1-z_2]\right),
\end{multline*}
the principle term in (\ref{eq:stat_phase}) with phase replaced with its first-order expansion. 
(If $\omega\in T^*\HH$ is a covector then by $\omega[z_1-z_2]$ we denote application of $\omega$ to vector $z_1-z_2$.) Let's prove it in a bit more details. We are going to deal with the $1/\tau$ order term in expansion in (\ref{eq:stat_phase}), the others are treated similarly. Put 
\begin{multline*}
	L_{\mathrm{PDO}}(z_1, z_2) = \tau^2\int_{\R} d\xi_1 \, \int_{\R} d\xi_2\, g(\Eta_{z_2,t}(\xi_1, \xi_2))\cdot \mathscr b_2(z_2, z_2,t, \Eta_{z_2,t}(\xi_1, \xi_2), \Theta_{z_2,t}(\xi_1, \xi_2))\times
	\\
	\times\left|\det\dfrac{\partial(\Theta_{z_2,t}, \Eta_{z_2,t})}{\partial(\xi_1, \xi_2)}\right|
	a(Q(z_2,t,\Eta_{z_2,t}(\xi_1, \xi_2),\Theta_{z_2,t}(\xi_1, \xi_2)))
	e^{i\tau (\xi_1\,dx+ \xi_2\,dy)[z_1-z_2]}. 
\end{multline*}
We took main term in integral for $L_t(\cdot,\cdot)$, replaced phase by Taylor expansion, put $z_1=z_2$ in both $a_1$ and $\mathscr b_2$, changed integration variables  as $(\theta,\eta)\mapsto \mathcal T_{z_2,t}(\theta,\eta)=(\xi_1,\xi_2)$  and, finally, multiplied by $\tau^3$. Our goal is to show that this constant-scale semiclassical ($\hbar=1/\tau$) PDO is indeed $\tau^3$ times principal term in $L_t(\cdot,\cdot)$. That is, we are going to show that if $$\tilde b(z_2,\xi_1,\xi_2) := (2\pi)^2\, \mathscr b_2(z_2, z_2,t, \Eta_{z_2,t}(\xi_1, \xi_2), \Theta_{z_2,t}(\xi_1, \xi_2))\cdot\left|\det\dfrac{\partial(\Theta_{z_2,t}, \Eta_{z_2,t})}{\partial(\xi_1, \xi_2)}\right|$$ then $b_{1,t}(z,\xi_1, \xi_2) := \tilde b(z, -\xi_1, -\xi_2)$ satisfies the requirements from the statement of our Proposition.

Denote $F := \Phi_\eta(z_1, f(z_0(z_1, z_2, t, \eta,\theta))) + \Phi_{(1)}(z_2, z_0(z_1, z_2, t, \eta,\theta),\eta)-2\varphi(t,\eta,\theta)$, this is the exponent in (\ref{eq:stat_phase}). Notice that $\Re F\le 0$ due to \cite[Lemma X.2.5]{Treves2}.
Observe that $\det\dfrac{\partial}{\partial(\theta,\eta)}\left.d_{z_1}F\right|_{z_1=z_2} \neq 0$. More carefully,  consider matrix
$$
\begin{pmatrix}
	\dfrac{\partial^2 F}{\partial \Re z_1 \, \partial \theta} & \dfrac{\partial^2 F}{\partial \Re z_1 \, \partial \eta} \\[4mm] 
	\dfrac{\partial^2 F}{\partial \Im z_1 \, \partial \theta} & \dfrac{\partial^2 F}{\partial \Im z_1 \, \partial \eta} 
\end{pmatrix}.
$$
If we write $\mathcal T_{z_2, t}(\theta,\eta)=\mathcal T_{z_2, t}^{(1)}\,dx+ \mathcal T_{z_2, t}^{(2)}\,dy\in T_{z_2}\HH$ then, by expression for (\ref{eq:diff_phase1}), the latter matrix, at $z_1=z_2$, is 
$$i
\begin{pmatrix}
	\dfrac{\partial \mathcal T_{z_2, t}^{(1)}}{\partial \theta} & \dfrac{\partial \mathcal T_{z_2, t}^{(1)}}{\partial \eta} \\[3mm] 
	\dfrac{\partial \mathcal T_{z_2, t}^{(2)}}{\partial \theta} & \dfrac{\partial \mathcal T_{z_2, t}^{(2)}}{\partial \eta} 
\end{pmatrix}
$$
which, by Lemma \ref{lemma:diffeo}, is non-degenerate.
Also $F|_{z_1=z_2}=0$ and therefore $d_{\theta,\eta}F|_{z_1=z_2}=0$. Then, for $z_1$ close to $z_2$, $|d_{\theta,\eta}F|\ge \const\cdot |z_1-z_2|$. 

Consider 
$$
I_1:=\int_0^{2\pi}d\theta\, \int_{\R^+} d\eta\, g(\eta) a_1(z_0(z_1, z_2,t,\eta,\theta)) \mathscr b_2(z_1, z_2,t,\eta,\theta) e^{\tau F},
$$
the original main term in  (\ref{eq:stat_phase}). Suppose that $|z_1-z_2| \ge \tau^{-2/3}$ but $|z_1-z_2|$ is small enough (less than some constant not depending on $\tau$). Let us show that $I_1$ is small then. In Euclidean coordinates $(\Re z_{1,2}, \Im z_{1,2})$, assume that $z_1-z_2$ lies in some cone thin enough. 
We may localize $I_1$ by multiplying its amplitude by a partition of unity. This allows us to assume that there exists a unit vector $v$ in $(\theta,\eta)$-plane such that 
\begin{equation}
	\label{eq:diff_sub}
	|d_{\theta,\eta}F[v]| \ge \const\cdot |z_1-z_2|
\end{equation}
on the whole support of integrand. In localization of $I_1$ change variables such that this integral will be sliced in $v$-direction. Now we apply repeated integration by parts in {$v\mbox{-direction}$} to show that $I_1=O(\tau^{-N})$ for any $N$. 

More formally, put $z=z_1-z_2=r e^{i\phi}$, $G(z):= \dfrac{\partial F}{\partial v}$ and integrate localized $I_1$ by parts in $v$-direction. We have $\tau G$ in denominator after this. Since $G|_{z=0}=0$, we have $G(z) = r\int_0^1 \langle\nabla_zG({\rho z}), e^{i\phi}\rangle\, d\rho$. Thus, $G/r$ is smooth in $v$-direction and, by (\ref{eq:diff_sub}), is separated from zero with $z$ small. This allows further integration by parts and leads to estimate $I_1=O(\tau^{-N})$ for any $N$. If we  assume that $z_1, z_2$ range a compact set which is not far from $\supp a$, and $|z_1-z_2|\ge \tau^{-2/3}$ then we have $O(\tau^{-N})$ estimate for the amount of such $z_1$, $z_2$ to quadratic form at the left-hand side of (\ref{eq:quadratic_form_equiv}).

The same concerns $L_{\mathrm{PDO}}$ (just integrate it by parts in the appropriate direction in $(\xi_1, \xi_2)$-plane).

Now assume that $|z_1-z_2| \le \tau^{-2/3}$. In $I_1$, apply Taylor expansion by degrees of $\Re(z_1-z_2)$ and $\Im(z_1-z_2)$ in $a_1(z_0(z_1, z_2,t,\eta,\theta))$ and in $\mathscr b_2$. Also, write $F$ by Taylor at $z_1=z_2$ and write 
$$
e^{\tau F} = e^{\tau\cdot i\mathcal T_{z_2,t}(\theta,\eta)[z_1-z_2]}\cdot e^{\tau \cdot O(|z_1-z_2|^2)}.
$$
Write long enough expansion 
for the remainder $\tau \cdot O(|z_1-z_2|^2)$. Further,
expand the second factor $e^{\tau \cdot O(|z_1-z_2|^2)}$ by Taylor (the exponent  is $o(1)$ therein). Degree of the expansions above can be taken large enough such that all the remainders are $O(1/\tau^4)$ which fits into precision claimed in the statement of our Proposition. We thus obtain 
an asymptotic expression which is  sum of terms like
$$
1/\tau\cdot (\Re(z_1-z_2))^\alpha (\Im(z_1-z_2))^\beta \tau^\gamma e^{\tau\cdot i\mathcal T_{z_2,t}(\theta,\eta)[z_1-z_2]}\times (\mbox{some smooth amplitude at }z_2).
$$
In all terms except for the main one 
we have $\alpha+\beta > \gamma$.  In $I_1$, change variable as $(\theta,\eta)\mapsto (\xi_1, \xi_2) = \mathcal T_{z_2,t}(\theta,\eta)$. In all terms in expansion 
except the main one we may integrate by parts $\alpha$ times with respect to $\xi_1$ and $\beta$ times with respect to $\xi_2$. Calderon--Vailliancourt Theorem then implies that all summands except for the main one bring to the original  operator terms whose $\|\cdot\|_{L^2\to L^2}$-norms are $O(\tau^{-4})$. The principal term leads to the proposed asymptotics.

The latter arguments also concern the case when $z_1$ is close enough to $z_2$ and both range a compact set. 
We also need to show that $L_t(z_1, z_2)$ gives a small operator in $L^2\to L^2$ when $z_1$ and $z_2$ are separated one from another by a positive constant or when at least one of them is far from $f^{-1}(\supp a)$. These cases are not covered by the above arguments. It is enough to estimate each $L_{t,\eta,\theta}(z_1, z_2)$. 

If at least one of $z_1$ or $z_2$ is far away from $\supp a$ then we make use of condition $\sup\limits_{n\in\mathbb N, \, z\in \HH} \|u_n\|_{L^1(\mathcal B_\HH(z,1))}<+\infty$. We may apply the similar uniform estimate on any \emph{Carleson square} which is $[x,x+y]\times[y,2y]\subset \HH$ for some $x\in\R$, $y>0$. Let $\supp a$ be close to, say, $h_{-it}(i,\pi)$.  
We need to show that, for any $C>0$, 
$$
\int\limits_{\{\dist(z_1,i)>R\}}
L_t(z_1,z_2)u(z_1)\bar u(z_2)\,d\mathcal A_2(z_1) \,  d\mathcal A_2(z_2) 
$$
can be forced to be $ O(\exp(-C\tau))$ by an appropriate choice of $R\gg 1$; and we need the same for ${\{\dist(z_1,i)>R\}}$ replaced with ${\{\dist(z_2,i)>R\}}$. 
%
%
To this end we apply argument similar to that of the technical part in Lemma~\ref{lemma:SelbergTransform}. It is useful  to notice that 
$\left|\dfrac{z_1-\tilde Z(P)}{\bar z_1 - Z(P)}\right|$ is bounded from the above and separated from zero for $P\in\supp a$, this is achieved by making $\supp a$ close enough to $h_{-it}(i,\pi)$. Then argue in a manner generally similar to the 
proof of the first assertion of Lemma~\ref{lemma:SelbergTransform}. 
We cut $\HH\setminus ([-C_1,C_1]\times[1/C_1,C_1])$ with $C_1\gg 1$ into Carleson squares, then estimate maximum of $\exp(-\tau c_t \cosh\dist(z, P))$ when $z$ ranges any of over each of the squares. 
We omit this technicality in our exposition.


If $z_1$ and $z_2$ are not far from  $\supp a$ but separated then we integrate by a bounded set when evaluating the operator. From the proof of the second assertion of Lemma~\ref{lemma:kernel_define_and_max} we see that either $\left(\dfrac{z_1-\tilde Z(P)}{\bar z_1 - Z(P)}\right)^\tau e^{-\tau c_\eta \cdot \cosh \dist(z_1,P)-\tau \varphi(t,\eta,\theta)}$ is small or the same for $z_2$ (absolute value is $\le e^{-c\tau}$ with some $c$ positive). 
This gives the desired.

Now, let's check the signs. Form given by PDO some symbol $\mathscr s$ is
\begin{multline*}	\dfrac{\tau^2}{(2\pi)^2}\int\limits_{\HH}d\mathcal A_2(z_1) \, \int\limits_{\HH} d\mathcal A_2(z_2) \, u(z_1) \bar{u}(z_2) \int\limits_{\R} d\xi_1' \, \int\limits_{\R} d\xi_2' \, \mathscr s(z_2, \xi_1', \xi_2') e^{i\tau (\xi_1'dx+\xi_2'dy)[z_2-z_1]}
=\\=
\dfrac{\tau^2}{(2\pi)^2}\int\limits_{\HH}d\mathcal A_2(z_1) \,\int\limits_{\HH}d\mathcal A_2(z_2) \, u(z_1) \bar{u}(z_2) \int\limits_{\R} d\xi_1 \, \int\limits_{\R} d\xi_2 \,\mathscr s(z_2, -\xi_1, -\xi_2) e^{i\tau (\xi_1dx+\xi_2dy)[z_1-z_2]}.
\end{multline*}
If $\mathscr s$ is as in the statement of our Proposition then the latter is $\int\limits_{\HH}d\mathcal A_2(z_1) \, \int\limits_{\HH} d\mathcal A_2(z_2) \, u(z_1) \bar{u}(z_2) L_{\mathrm{PDO}}(z_1, z_2)$. But $\tau^{-3}L_{\mathrm{PDO}}$ gives the main part of kernel $L_t(z_1, z_2)$. Proof is complete.
$\blacksquare$

\end{document}